\documentclass[final,sort&compress,5p,times]{elsarticle}

\usepackage{microtype}
\usepackage{flushend}

\usepackage[cmex10]{amsmath}
\usepackage{amsthm}
\usepackage{amsfonts}
\usepackage{wasysym}
\usepackage{mathrsfs}
\usepackage{graphicx}
\usepackage{float}
\usepackage{tikz}
\usetikzlibrary{shapes.geometric}

\newtheorem{theorem}{Theorem}

\newtheorem{lemma}{Lemma}

\theoremstyle{definition}

\newtheorem{assumption}{Assumption}

\theoremstyle{remark}
\newtheorem{remark}{Remark}

\journal{arXiv}

\begin{document}

\begin{frontmatter}

 \title{Control of Positive Systems with an Unknown State-Dependent Power Law Input Delay and Input Saturation\tnoteref{t1}}
\author[address1]{Damon E.~Ghetmiri\fnref{fn1}}
\fntext[fn1]{\emph{Now at:} ASML, 17075 Thornmint Court, San Diego, CA 92127-2413, USA}

\author[address1]{Amor A.~Menezes\corref{cor1}}\ead{amormenezes@ufl.edu}
\tnotetext[t1]{This paper was not presented at any IFAC meeting.}
\cortext[cor1]{Corresponding author.}

\address[address1]{Department of Mechanical and Aerospace Engineering, University of Florida, Gainesville, FL 32611-6250, USA}                  

\begin{abstract}                  
This paper is motivated by a class of positive systems with an input that is subject to an unknown state-dependent power law delay as well as saturation. For example, biological networks have non-negative protein concentration states. Mass action kinetics in these systems result in power law behavior, while complex interactions cause signal propagation delays. Incomplete network characterization makes delay state-dependence unknown. Manipulating network activity via modulated protein concentrations to attain desired performance is restricted by upper-bounds on concentration actuator authority. Here, an innovative control law exploits system dynamics to compensate for control domain restrictions. A Lyapunov stability analysis establishes that the reference tracking error of the closed-loop system is uniformly ultimately bounded. Numerical simulations on a human coagulation model show controller efficacy and better performance compared to the relevant literature. This example application steps toward personalized, closed-loop treatments for trauma coagulopathy, which currently has 30\% mortality with open-loop clinical approaches.
\end{abstract}


\end{frontmatter}


\section{Introduction}
Positive systems \cite{farina2011positive,kaczorek2012positive,blanchini2015switched,rantzer2018tutorial,rantzer2021scalable} are a class of dynamical systems whose states are confined to the positive orthant of $\mathbb{R}^n$. These systems model many real-world applications, including communication networks \cite{shorten2006positive}, ecosystems \cite{murray2002mathematical}, thermodynamics \cite{haddad2010nonnegative}, biological systems \cite{blanchini2014piecewise,briat2016antithetic}, and physiological systems \cite{hovorka2004nonlinear}. Intrinsic to these systems are non-negative variables like information, population levels, absolute temperature, and substance concentration. Since both system states and control inputs are restricted to the positive orthant, many classical control ideas are not applicable \cite{shu2008positive} or not easily adapted \cite{shen2016static}. Nevertheless, some controllability and reachability results exist   \cite{coxson1987positive,benvenuti2004tutorial,valcher2009reachability,guiver2014positive,eden2016positive,briat2020biology}. 

Control systems can have time delays \cite{wu2010stability} that are nonlinear \cite{otto2019nonlinear}. Dynamical systems with input time delays may exhibit oscillations, instability, and poor performance \cite{hale1993introduction}. Techniques to control systems with input delays \cite{park2014stability} have been developed for delays that are constant \cite{li1997delay,sheng2021switched}, time-varying \cite{obuz2017unknown}, and state-dependent \cite{bekiaris2012compensation}. Systems with an unknown time delay have also been studied. The delays for these systems are constant \cite{alibeji2017pid,deng2021state}, time-dependent and bounded by known constants \cite{koo2020output,nguyen2021state}, or of small magnitudes  \cite{cai2018adaptive}. Moreover, practical control systems also have actuator nonlinearities \cite{zhang2015exact}. Saturation constraints on the magnitude of the input control signal can severely limit control performance and lead to system instability \cite{perez2009saturation}.

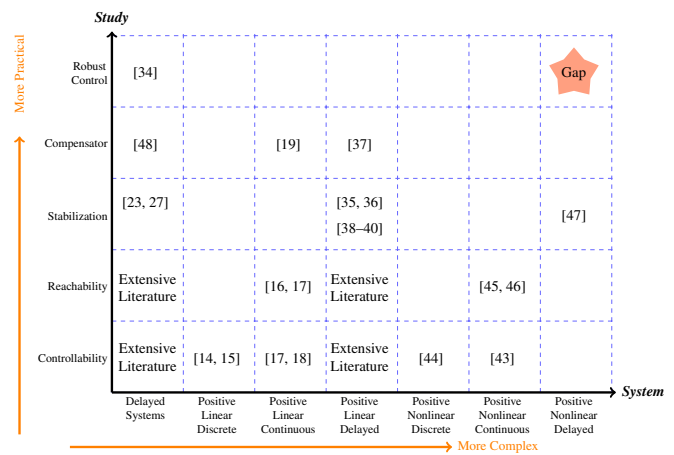
\begin{figure}[h!]
    \centering
    \resizebox{0.483\textwidth}{!}{
    \begin{tikzpicture}[every text node part/.style={align=center}]
\draw[help lines,step=1.7cm, color=blue!60, dashed] (0,0) grid (12,8.6);
\draw[->,ultra thick] (0,0)--(12,0) node[right]{\textbf{\textit{System}}};
\draw[->,ultra thick] (0,0)--(0,8.6) node[above]{\textbf{\textit{Study}}};

\small
\draw[->,color=orange,ultra thick] (-1,-1.3)--(8.1,-1.3)  node[color=orange,right]{More Complex};
\draw[->,color=orange,ultra thick] (-2.2,-1)--(-2.2,6.1);
\node[color=orange,rotate=90] at (-2.2,7.6) {More Practical};

\footnotesize
\draw (.8,0) node[below]{ Delayed \\ Systems};
\draw (2.5,0) node[below]{ Positive \\ Linear \\Discrete};
\draw (4.2,0) node[below]{ Positive \\ Linear \\ Continuous};
\draw (5.9,0) node[below]{ Positive \\ Linear \\ Delayed};
\draw (7.6,0) node[below]{ Positive \\ Nonlinear \\ Discrete};
\draw (9.3,0) node[below]{ Positive \\ Nonlinear \\ Continuous};
\draw (11,0) node[below]{ Positive \\ Nonlinear \\ Delayed};
\normalsize

\footnotesize
\draw (0,.8) node[left]{Controllability};
\draw (0,2.5) node[left]{Reachability};
\draw (0,4.2) node[left]{Stabilization};
\draw (0,5.9) node[left]{Compensator};
\draw (0,7.6) node[left]{Robust \\ Control};
\normalsize

\draw (0.8,4.5) node{\cite{park2014stability,bekiaris2012compensation}};
\draw (0.8,5.9) node{\cite{sabatier2020power}};
\draw (0.8,7.6) node{\cite{perez2009saturation}};

\draw (2.5,0.8) node{\cite{coxson1987positive,benvenuti2004tutorial}};

\draw (4.2,0.8) node{\cite{guiver2014positive,eden2016positive}};
\draw (4.2,2.5) node{\cite{valcher2009reachability,guiver2014positive}};
\draw (4.2,5.9) node{\cite{briat2020biology}};


\draw (5.9,4.5) node{\cite{rami2007positive,liu2010stability}}; 
\draw (5.9,3.9) node{\cite{zhu2013exponential,zhang2015h,hong2019optimization}};
\draw (5.9,5.9) node{\cite{liu2009constrained}};

\draw (7.6,0.8) node{\cite{naim2018controllability}};

\draw (9.3,0.8) node{\cite{klamka1996constrained}};
\draw (9.3,2.5) node{\cite{greenstreet1999reachability,asarin2003reachability}};

\draw (11,4.2) node{\cite{liu2015stability}};

\draw (.85,.85) node{Extensive \\ Literature};
\draw (5.9,2.5) node{Extensive \\ Literature};
\draw (.85,2.5) node{Extensive \\ Literature};
\draw (5.9,.85) node{Extensive \\ Literature};

\node[shape=star,fill=orange!50!red!50!white,inner sep=0pt] at (11,7.6) {Gap};

\end{tikzpicture}
}
\vspace{-2em}
    \caption{White space chart categorizing the relevant positive systems literature.
    }
    \label{fig:whitespacechart}

\end{figure}

For positive linear time delay systems, their stability \cite{rami2007positive,liu2010stability}, controller design \cite{liu2009constrained}, and output feedback stabilization \cite{zhu2013exponential,zhang2015h,hong2019optimization} have been studied, but the delays in such systems are constant \cite{rami2007positive}, known \emph{a priori} \cite{arogbonlo2019design}, or have a small known time-varying bound (and are hence of small magnitude) \cite{liu2010stability, zhu2013exponential,hong2019optimization}. In real world applications, the delay magnitude is often unavailable \emph{a priori}, and is particularly problematic when it is of the same order as the system's time-constant \cite{ajmeri2015simple}. For positive nonlinear systems, results on controllability \cite{klamka1996constrained,naim2018controllability}, reachability \cite{greenstreet1999reachability,asarin2003reachability}, and stabilization with delays \cite{liu2015stability} exist. Fig.~\ref{fig:whitespacechart} classifies the relevant positive systems literature according to complexity and practicality, and shows where the current research gap is.

\vspace{-0.25em}

A common nonlinear time delay, $\tau$, is the power law ($\tau(x)=\gamma x^{-k}$ for some constants $\gamma$ and $k$), which affects battery electrochemistry, thermal conduction, wave propagation in viscoelastic materials, acoustics, electrical distribution networks, and biochemical pathways \cite{sabatier2020power}. There exist positive systems with inputs that are subject to nonlinear power law time delays and saturation. In biochemical pathways, non-negative protein concentrations are modulated in a chemical interaction network, and their effects are propagated downstream depending on intermediate network activity. Innate mass action kinetics cause delays that can be captured by power laws \cite{sun2004explanation}, and activity saturation can be captured by a logistic function \cite{meyer1994bi}. Power law time delays are also seen at the cellular level: after being confronted with antibiotics, \emph{Escherichia coli} cells rejuvenate with a power law time delay that stems from bacterial persistence \cite{csimcsek2019power}. 

This paper describes how to close the loop for positive systems with simultaneous input saturation and a nonlinear power law state-dependent input delay, a combination that has not been previously studied. The time delay magnitude can be the same order as the system time-constant, and the delay parameters can be unknown. Intellectual contributions include:
\begin{enumerate}[(1)]
    \item a control law for positive systems with a state-dependent power law input delay and input saturation that takes advantage of natural dynamics to compensate for restrictions on the domain of feasible control; 
    \item a Lyapunov analysis that uses Lyapunov-Krasovskii functionals to establish that the reference tracking error of the resulting closed-loop system is uniformly ultimately bounded;
    \item a controller for a model of the human coagulation positive system, to meet the need for one that can remedy coagulation deficits in trauma patients despite a state-dependent power law input delay and input saturation.
\end{enumerate}

In what follows, Section \ref{Section:ProblMotiv} provides a real-world motivating application. Section \ref{section:Theoretical} develops a controller that satisfies system constraints. Section \ref{section5} confirms suitable controller performance on the application model. Section \ref{section6} states conclusions.

\section{Biomolecular Problem Motivation} \label{Section:ProblMotiv}
A biomedical example with an input power law time delay and logistic saturation is the human coagulation system, for which one possible state-space realization from the third-order transfer function model in the literature \cite{menezes2017targeted} is:
\begin{equation} \label{coagulationmodel}
    \begin{split}
        &\dot{x}_1(t)=x_2(t)-d_1 x_1(t), \\
        &\dot{x}_2(t)=x_3(t)-d_2 x_2(t), \\
        &\dot{x}_3(t)=-d_3x_3(t)+bv(t-\tau).
    \end{split}
\end{equation}
In \eqref{coagulationmodel}, state $x_1(t)$ is the concentration of thrombin (factor IIa), a key end product protein of coagulation that is both anticoagulant and procoagulant \cite{narayanan1999multifunctional}. State $x_2(t)$ is the concentration of protein complex prothrombinase, which includes factor Xa. State $x_3(t)$ is the concentration of a protein complex of tissue factor (TF) bound to another protein, factor VIIa. The terms with scalar parameters, $d_1>0$, $d_2>0$, and $d_3>0$, represent degradation of the respective proteins. Input parameter $b$ is also a scalar.

TF initiates coagulation after perforating vascular injury. Therefore, TF is an input, $v(t)$, in \eqref{coagulationmodel}. TF is a membrane protein that, when exposed, binds with blood plasma factor VIIa, as represented by state $x_3(t)$. The bound-protein complex TF-VIIa is a key driver of coagulation through activating factor IX into IXa and factor X into Xa. Hence, coagulation continues as long as TF is released in the blood. Output thrombin can be measured by the Calibrated Automated Thrombogram (CAT).

The human coagulation system is a positive system. Consider the linear time-invariant form of \eqref{coagulationmodel}, $\dot{x}(t)=Ax(t)+Bv(t-\tau)$. The resultant matrix $A$ is a Metzler matrix, which is a square real matrix $\mathbb{M}$ whose off-diagonal elements are nonnegative, i.e., $A=(a_{ij});~a_{ij}\geq 0 ~\forall~ i\not=j$ \cite{luenberger1979introduction}.
A necessary and sufficient condition for a linear state space model to be positive is for $A$ to be a Metzler matrix \cite{kaczorek2012positive,rantzer2018tutorial}, and so \eqref{coagulationmodel} is a positive system. 

\begin{figure}[h!]
    \centering
    \includegraphics[width=.47\textwidth]{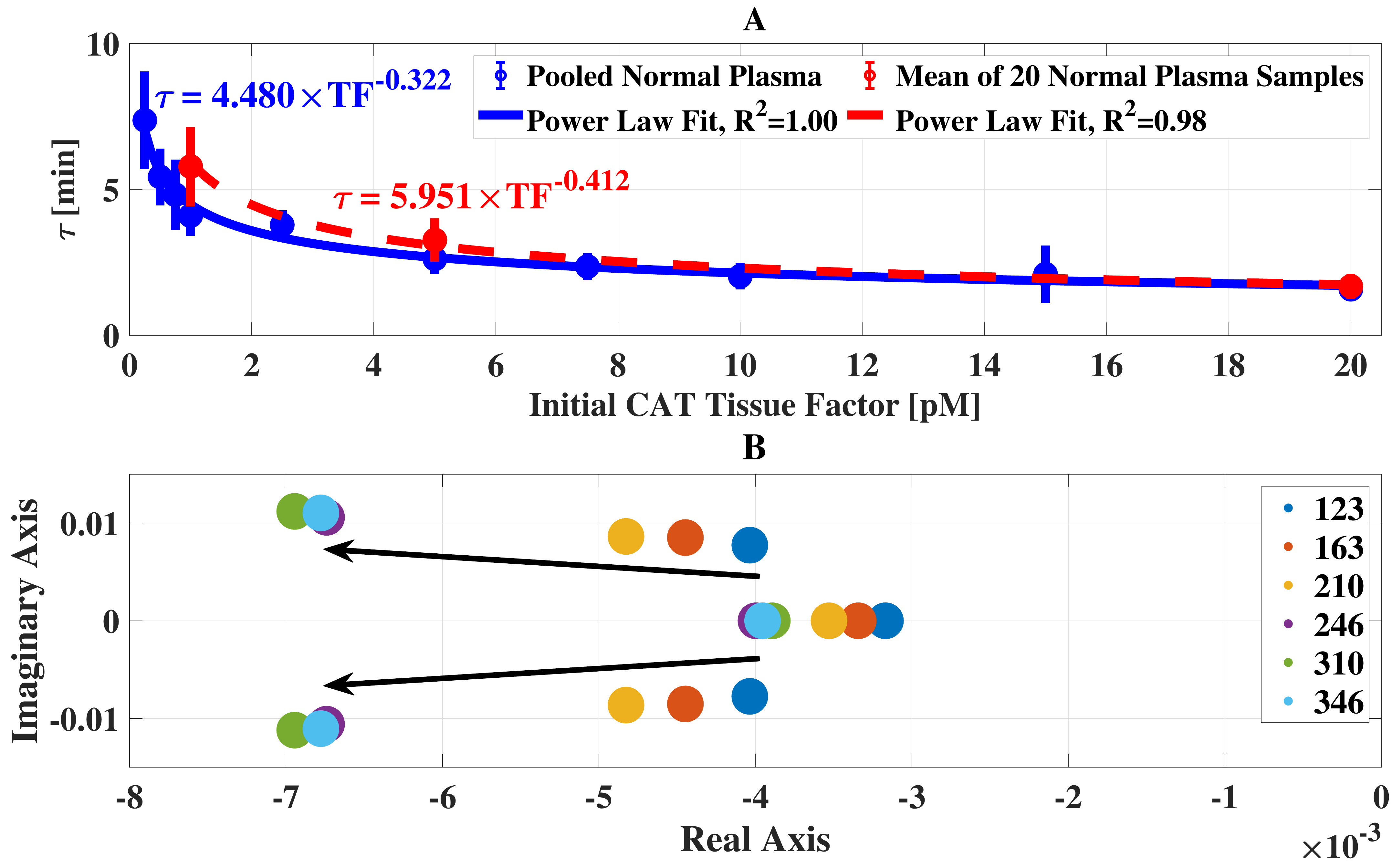}
    \vspace{-1em}
    \caption{Time delay and saturation in the human coagulation model \eqref{coagulationmodel}. (A) Time delays extracted from 
    \cite{menezes2017targeted} 
    are correlated with TF concentration via a power law, $\tau=\gamma\times TF^{-k_d}$. Red indicates biological replicates from 20 different normal plasma samples ($R^2=1.00$), blue corresponds to technical replicates on a pool of normal plasma samples ($R^2=0.98$). (B) Increases in the concentration of factor X (a coagulation protein) move the poles of a linear dynamical model of the coagulation system 
    to a certain limit. Beyond that, concentration increases have little 
    effect, which confirms protein input saturation. Similar pole-movement and saturating behavior exists for other coagulation proteins \cite{ghetmiri2021personalized}.}
    \label{fig1}
\end{figure}

Experimental results \cite{menezes2017targeted} show that TF causes clotting that is subject to a nonlinear delay, $\tau$. Fig.~\ref{fig1}A plots $\tau$ and a power law fit for plasma that was pooled from normal human samples, and also $\tau$ and a power law fit for the mean of 20 different normal human plasma samples that were assessed individually. The power law fits are a function of TF, 
and hence 
of state $x_3$ in \eqref{coagulationmodel}. 

Clotting activity can be modulated by adding coagulation proteins (e.g., TF, factor II, factor VIII, and factor X) to increase their concentrations. Fig.~\ref{fig1}B (data from \cite{menezes2017targeted}) shows how the system poles of \eqref{coagulationmodel} move upon the experimental addition of concentrations of factor X. The plot indicates that adding coagulation proteins has diminishing returns due to protein activity saturation. The complex pole-pair location in Fig.~\ref{fig1}B 
stops moving left despite the continued addition of factor X \cite{ghetmiri2021personalized}. 

Closed-loop addition of protein concentrations to rectify abnormal clotting, despite input time delay and saturation, can personalize treatment for trauma patients \cite{gonzalez2016trauma}, 
and is considered
beneficial \cite{ghetmiri2021personalized}. 
The control objective here is for thrombin to follow a desired reference trajectory to achieve normal clotting following trauma \cite{ghetmiri2021personalized, ghetmiri2021controltheoretic,ghetmiri2022nonlinear}. This is an important endeavor because trauma is the leading cause of death and disability in the United States between ages 1 and 44 \cite{sleet2011injury}, with bleeding 
a major cause of these deaths \cite{moore2021trauma}.The mortality from massive transfusion remains high at 30\% \cite{cantle2017prediction}. 

\section{Theoretical Analysis}\label{section:Theoretical}

\subsection{Problem Formulation}\label{subsec:ProblemForm}
Let $\mathbb{R}$ denote the set of real numbers; $\mathbb{R}^n$ and $\mathbb{R}^{n \times m}$ denote the $n$-dimensional and $n\times m$-dimensional real spaces, respectively; and $\mathbb{R}_+$ denote the set of non-negative real numbers. 
%
We note that \eqref{coagulationmodel} is a cascading system \cite{allaby2013dictionary}, a type of dynamical system that is characterized by the transfer of mass and energy along a chain of component sub-systems, i.e., the output from one sub-system becomes the input for an adjacent sub-system. These systems have broad applicability in biology \cite{young2017dynamics}, power systems \cite{guo2017critical}, and geomorphology \cite{allaby2013dictionary}. To generalize \eqref{coagulationmodel}, we consider the $n^{th}$-order cascading dynamical system: 
\begin{align} 
        \dot{x}_i(t) &= x_{i+1}(t)-d_ix_i(t)+\mathcal{H}_i(\mathbf{x}(t)),~~~~~~i = 1,2,\ldots,n-1;\nonumber\\ 
        \dot{x}_n(t) &=-d_nx_n(t)+\mathcal{F}\left(\mathbf{x}\left(t\right)\right)+g\left(u\left(t-\tau\left(x_n\left(t\right)\right)\right)\right)\label{dymodel}, 
\end{align}
where $\mathbf{x}(t)\in\mathbb{R}_+^n$ is an $n$-vector of non-negative system states; $\mathbf{x}(0)=\mathbf{x_0}$ is the initial condition; $x_i(t)$ is the $i^{th}$ system state; $d_i\in\mathbb{R}_+$ is a non-negative constant for each $i^{th}$ state; $\mathcal{H}_i:\mathbb{R}_+^n\rightarrow\mathbb{R}_+$ is a non-negative map for each $i^{th}$ state; $\mathcal{F}:\mathbb{R}_+^n\rightarrow\mathbb{R}_+$ describes the dynamics of the $n^{th}$ state;  $\tau:\mathbb{R}_+\rightarrow\mathbb{R}_+$ is a state-dependent, time-varying, unknown, non-negative input delay; $u:\mathbb{R}_+\rightarrow\mathbb{R}_+$ is a non-negative system input; and $g(u):\mathbb{R}_+\rightarrow\mathbb{R}_+$ is a non-negative saturation of the input $u$ that we model as a logistic function, equivalently a hyperbolic tangent function:
\begin{equation} \label{saturationFuncg}
    g(u)=\frac{\beta}{1+e^{-k_s\left(u-\eta\right)}}=\frac{\beta}{2}+\frac{\beta}{2}\tanh{\left(\frac{k_s}{2}\left(u-\eta\right)\right)},
\end{equation}
where $\beta, k_s,$ and $\eta$ are positive constants, and $\eta \geq \beta$. We can write the dynamical system \eqref{dymodel} as:\small
\begin{equation}
    \dot{\mathbf{x}}=\underbrace{\begin{bmatrix}
    -d_1 &1 &0 &\ldots &0\\
    0 &-d_2 &1 &\ldots &0\\
    \ldots\\
    0 &\ldots &0 &-d_{n-1} &1\\
    0 &&\ldots &0&-d_n
    \end{bmatrix}}_A
    \mathbf{x} + \begin{bmatrix}
    \mathcal{H}_1(\mathbf{x}) \\ \mathcal{H}_2(\mathbf{x})\\ \vdots \\ \mathcal{H}_{(n-1)}(\mathbf{x}) \\ \mathcal{F}(\mathbf{x}) 
    \end{bmatrix}
    + \begin{bmatrix}
    0_{(n-1) \times 1} \\ 1
    \end{bmatrix} g(u). \label{eq:extendedMetzler}
\end{equation}\normalsize 
Here, $A\in \mathbb{M}_n$ is a Metzler matrix. 
%
For our generalization \eqref{dymodel} to be a positive system, we use the known result \cite{kaczorek2015analysis} that 
a continuous-time nonlinear system
\begin{equation*}
        \dot{\mathbf{x}} = A\mathbf{x}+q(\mathbf{x},\mathbf{u}),
\end{equation*}
is positive if and only if 
\begin{equation}
    A\in \mathbb{M}_n \text{ and } q(\mathbf{x},\mathbf{u})\in\mathbb{R}_+^n, \forall \mathbf{x}\in\mathbb{R}_+^n,~ \mathbf{u}\in\mathbb{R}_+^n,~ t\geq 0, \label{eq:nonlinear_positive}
\end{equation}
%
%
which is true by definition:
$\mathcal{H}_{i}(\mathbf{x}(t))\in \mathbb{R}_+$, $i = 1,2,\ldots,n-1$, $  \mathcal{F}(\mathbf{x}(t))\in\mathbb{R}_+$, and $g(u)\in\mathbb{R}_+$, $\forall x\in\mathbb{R}_+^n,~ u\in\mathbb{R}_+^n,~ t\geq 0$. We can confirm that the input to \eqref{dymodel} is non-negative, i.e., in \eqref{saturationFuncg} we have $g(u) \geq 0$. This is because $-1 \leq \tanh(\cdot) \leq 1 \Rightarrow -\frac{\beta}{2}\leq\frac{\beta}{2}\tanh(\cdot)\leq\frac{\beta}{2} \Rightarrow 0 \leq \frac{\beta}{2}+\frac{\beta}{2}\tanh(\cdot) \leq \beta$. Non-negativity will also be enforced in controller design. 
Therefore, 
for any initial $\mathbf{x_0} \in \mathbb{R}_+^n$ and 
any $g(u)$, 
we have that $\mathbf{x}(t)\in \mathbb{R}_+^n$ $\forall t\geq 0$. 

The input in \eqref{dymodel} has a state-dependent time delay $\tau$ that is an unknown power law
\begin{equation}
    \tau(x_n)=\gamma x_n^{-k_d}, \label{eq:powerlaw}
\end{equation}
where $\gamma$ and $k_d$ are positive constants (e.g., Fig.~\ref{fig1}A). 
We take $\tau$ to be a function of $x_n$ alone (e.g., Fig.~\ref{fig1}A). We use the notation 
\begin{equation*}
u_{\tau}\triangleq\begin{cases}
u(t-\tau), &\text{ if } (t-\tau)>0,\\
0, &\text{ otherwise.}
\end{cases}
\end{equation*}
\normalsize

As state $x_n(t)$ goes to zero, the time delay increases. Given $x_n(t)\geq0$ as \eqref{dymodel} is a positive system, to keep $\tau$ finite:
\begin{assumption} \label{assumption delay infinite state}
We assume that state $x_n(t)$ is strictly positive, i.e., there exists an arbitrarily small $\varphi_1 \in \mathbb{R}_+$ such that $x_n(t) > \varphi_1$ $\forall t\geq0$.
\end{assumption}
\noindent When considering our motivating problem, this is a reasonable assumption since TF drives coagulation and must have a concentration that is greater than zero to do so.

Assumption \ref{assumption delay infinite state} only serves to keep $\tau$ finite, i.e., using \eqref{eq:powerlaw} we have that $\tau <  \gamma\varphi_1^{-k_d}$. With $\varphi_1$ arbitrarily small, $\gamma\varphi_1^{-k_d}$ is arbitrarily large. The state bound $\varphi_1$ in the assumption is unknown and unspecified. Thus, the time delay $\tau$ can also be unknown. But:

\noindent\begin{assumption} \label{delayassumption}
We assume that the input delay $\tau$: 
\begin{enumerate}[(a)]
\item is differentiable; and 
\item varies slow enough, i.e., $\exists~\varphi_2 \in \mathbb{R}_+$, so that $|\dot{\tau}|< \varphi_2$ $\forall t\geq 0$.
\end{enumerate}
\end{assumption}

We choose any $\varphi_2>d_nk_d\gamma\varphi_1^{-k_d}$. As above, if $\varphi_1$ is arbitrarily small, then $\varphi_2$ is arbitrarily large, and thus Assumption~\ref{delayassumption}(b) is not restrictive. Let $h_{i}(\mathbf{x}(t)) = -d_ix_i(t)+\mathcal{H}_i(\mathbf{x}(t))$, and let $f(\mathbf{x}(t)) = -d_nx_n(t)+\mathcal{F}(\mathbf{x}(t))$. To facilitate analyzing our nonlinear system, we make:
\noindent\begin{assumption} \label{assumption_FunctionBound}
We assume that the dynamics \eqref{dymodel} are such that: 
\begin{enumerate}[(a)]
\item functions $f$, $g$, and $h_i ~ (i=1,2,\ldots,n-1)$ are differentiable; 
\item the function $f$ and its first partial derivative is bounded; and
\item the functions $h_i$ and their first $n-i+1$ partial derivatives are bounded.
\end{enumerate}
\end{assumption}

Let $x_r(t)\in \mathbb{R}_+$ be a desired state trajectory that satisfies:
\begin{assumption} \label{referenceassumption}
We assume that the reference trajectory $x_r(t)$ is such that all its time  derivatives 
$\dot{x}_r(t),\ddot{x}_r(t),\ldots$ exist and are bounded by positive constants for all $t\geq 0$.
\end{assumption}
\noindent The control problem that we have to solve is then: Design $u(t)$ so that, for some $T\in\mathbb{R}_+$, $\exists\epsilon \in \mathbb{R}_+$: $|x_r(t)-x_1(t)|\leq\epsilon, \forall t\geq T$, i.e., state $x_1(t)$ of \eqref{dymodel} tracks $x_r(t)$ within $\epsilon$ for all $t\geq T$.

\subsection{Controller Development}\label{section3-2}
We 
can define a tracking error $e_1 (t)\in\mathbb{R}$ as
\begin{equation} \label{e1error}
    e_1(t) \triangleq x_r(t)-x_1(t).
\end{equation}
We can also define auxiliary tracking error signals \cite{xian2004continuous} $e_i(t)\in\mathbb{R}$, $i=2,3,\ldots,n$:
\small
\begin{align} \label{error sys en}
    {e}_2(t) &\triangleq \dot{e}_1(t)+e_1(t), \\
    {e}_3(t) &\triangleq \dot{e}_2(t)+e_2(t)+e_1(t), \\
    \ldots \nonumber \\
    {e}_n(t) &\triangleq \dot{e}_{n-1}(t)+e_{n-1}(t)+e_{n-2}(t). \label{enerror}
\end{align}
\normalsize
Let $x_i^{(j)}$ be the $j^{th}$ time derivative of $x_i(t)$. Then \eqref{dymodel} can be written 
\small
\begin{equation} \label{dymodelx1}
    x_1^{(n)}=f(\mathbf{x})+g(u_\tau)+\sum_{i=1}^{n-1}h_i^{(n-i)}(\mathbf{x}). 
\end{equation}
\normalsize
Similarly, an expression for $e_i(t)$ is
\small
\begin{equation} \label{general ei}
    e_i(t)=\sum_{j=0}^{i-1}a_{i,j}e_1^{(j)},
\end{equation}
\normalsize
where $a_{i,j}$ is \cite{xian2004continuous_incollection}:
\begin{align*}
    a_{i,0} = &\frac{1}{\sqrt{5}}\left(\left(\frac{1+\sqrt{5}}{2}\right)^i - \left(\frac{1-\sqrt{5}}{2}\right)^i \right),~~i=2,3,...,n, \\
    a_{i,j} = &\sum_{p=1}^{i-1} a_{i-p-j+1,0} a_{p+j-1,j-1},\\ &~~~~~~~~~~~~~~~~~~~~~i=3,4,...,n,~~j=1,2,...,(i-2),\\
    a_{i,i-1}=&1,~~~~~~~~~~~~~~~~~~i=1,2,...,n. 
\end{align*}

Defining the following error signal $e_u(t)\in \mathbb{R}$ will help obtain a delay-free input expression for the closed-loop error system,
\begin{equation} \label{eu error}
    e_u(t) \triangleq -\int_{t-\hat{\tau}}^t \dot{u}(\theta)d\theta. 
\end{equation}
This expression uses an estimate of time delay $\hat{\tau}\in\mathbb{R}_+$, but the quality of this estimate can be poor given a lack of knowledge about $\tau$. To generate $\hat{\tau}$, we can exploit the form of \eqref{eq:powerlaw}. We can take $\hat{\tau}=\hat{\gamma}(x_n(0))^{-\hat{k}_d}$, where $\hat{\gamma}$ and $\hat{k}_d$ are chosen positive constants, and $x_n(0)$ is the initial condition of the $n^{th}$ state. Thus, $\dot{\hat{\tau}}=0$. 
We define the quality of our estimate as $|\Tilde{\tau}| \leq \Bar{\Tilde{\tau}}$, where $ \Tilde{\tau}=\tau-\hat{\tau} $, and $\Bar{\Tilde{\tau}}\in\mathbb{R}_+$.

In the next subsection, after some preliminaries and using a Lyapunov stability analysis, we will show that a controller 
is 
\begin{equation} \label{controllereq}
u(t)\triangleq \Bigg(\operatorname{sgn}\bigg(\frac{\operatorname{sgn}\left(e_1\left(t\right)\right)+1}{2}\bigg)\Bigg) \Big(k\big( e_n(t) - e_n(t_0)\big)+\nu(t)\Big),
\end{equation}
where $\operatorname{sgn}(x)$
is the signum function, $e_n(t_0)\in\mathbb{R}$ is an initial error signal, $k\in\mathbb{R}_+$ is a designed positive constant, and $\nu(t)$ is the solution to the ordinary differential equation
\begin{equation} \label{controller nu}
    \dot{\nu}(t) = k(\lambda e_n(t) + \alpha e_u(t)),
\end{equation}
having initial condition $\nu(t_0)\in\mathbb{R}$, with $\lambda\in\mathbb{R}_+$ and $\alpha\in\mathbb{R}_+$ also being designed positive constants. 

We define another auxiliary error signal $e_a(t)\in\mathbb{R}$, with
\begin{equation} \label{e_a error}
    e_a(t) \triangleq \dot{e}_n(t)+\lambda e_n(t) + \alpha e_u(t).
\end{equation}
\begin{assumption} \label{assumption error bound}
We assume that the error system is bounded, i.e., $\|e_a(\mu)\|<\psi$, $\psi\in\mathbb{R}_+$. 
\end{assumption}
\noindent This follows since biological systems are globally stable, even if locally unstable at short time scales \cite{sole1992structural}. 

The function $\operatorname{sgn}(x)$ is continuous and differentiable everywhere except at the singular point $x=0$, but under the generalized notion of differentiation in distribution theory, the derivative of the signum function is $2\delta(x)$ \cite{bracewell1986fourier}, where $\delta(x)$ is the Dirac delta function. 
Hence, the derivative of \eqref{controllereq}, 
\small
\begin{equation} \label{udot}
\begin{split}
\dot{u}(t)=&\big(\delta\left(e_1\left(t\right)\right)\delta\left(\operatorname{sgn}\left(e_1\left(t\right)\right)+1\right)\dot{e}_1\left(t\right)\big)\big(k\left( e_n\left(t\right) - e_n\left(t_0\right)\right)+\nu\left(t\right)\big)\\ &+\operatorname{sgn}\left(\frac{\operatorname{sgn}(e_1\left(t\right))+1}{2}\right) k e_a\left(t\right) =\operatorname{sgn}\left(\frac{\operatorname{sgn}\left(e_1\left(t\right)\right)+1}{2}\right) k e_a\left(t\right),
\end{split}
\end{equation}
\normalsize
is defined everywhere. Note:  $\delta(e_1(t))\delta(\operatorname{sgn}(e_1(t))+1)=0~\forall e_1(t)$.


The novel structure of the control law \eqref{controllereq} uses the natural dynamics of the system through the signum function. Given \eqref{eq:extendedMetzler} and the necessary and sufficient condition of positive systems in \eqref{eq:nonlinear_positive}, this signum function ensures that the controller only boosts the system, $u\in\mathbb{R}_+^n$, when necessary (when $e_1 \geq 0$), which keeps the system states in the positive orthant while also ensuring tracking of the reference signal. Otherwise, if state reductions are required ($e_1 < 0$), the controller switches off and takes advantage of existing natural decay dynamics, since it cannot supply negative inputs by the definition of a positive system.

\subsection{Stability Analysis}\label{section4}
Here, we analyze the performance of the controller in Section \ref{section3-2}. 
Taking the derivative of \eqref{e_a error}, the dynamics for $e_a(t)$ are
\begin{equation}
    \dot{e}_a(t)=\ddot{e}_n(t)+\lambda \dot{e}_n(t) + \alpha \dot{e}_u(t).\label{eq:ea2}
\end{equation}
These dynamics can be obtained by substituting into \eqref{eq:ea2} the first time-derivatives of \eqref{enerror} and \eqref{dymodelx1}, the second time-derivative of \eqref{general ei} with $i=n$, and the $(n+1)^{th}$ time-derivative of \eqref{e1error}. In what follows, we simplify our notation by no longer indicating time dependence. Using $u_{\tau}\triangleq u(t-\tau(x(t))$, and computing
\begin{equation}\label{eq:doteu}
    \dot{e}_u = (1-\dot{\hat{\tau}})\dot{u}_{\hat{\tau}} - \dot{u}=\dot{u}_{\hat{\tau}} - \dot{u},
\end{equation}
we obtain: 
\begin{equation} \label{rdot1}
\begin{split}
\dot{e}_a=&\sum_{j=0}^{n-1}a_{n,j}e_1^{(j+2)}+\lambda\dot{e}_n+\alpha\dot{e}_u,\\
    =&\sum_{j=0}^{n-2}a_{n,j}e_1^{(j+2)}+x_r^{(n+1)}-\dot{f}(\mathbf{x})+\lambda\dot{e}_n-\sum_{i=1}^{n-1}h_i^{(n-i+1)}(\mathbf{x}) \\
    &+\left(\alpha-\frac{\beta k_s}{4}\left( 1-\tanh^2\left(\frac{k_s}{2}\left(u_{\tau}-\eta\right)\right) \right) (1-\dot{\tau})\right) \dot{u}_{\tau} \\ 
    & +\alpha\left( \dot{u}_{\hat{\tau}} - \dot{u}_{\tau} \right)- \alpha\dot{u}.
\end{split}
\end{equation}
Substituting \eqref{udot} into \eqref{rdot1} we obtain

\footnotesize
\begin{equation} \label{rdot2}
    \begin{split}
        \dot{e}_a&=\sum_{j=0}^{n-2}a_{n,j}e_1^{(j+2)}+x_r^{(n+1)}-\dot{f}(\mathbf{x})-\sum_{i=1}^{n-1}h_i^{(n-i+1)}(\mathbf{x})+\lambda\dot{e}_n \\
        &+\left(\alpha-\frac{\beta k_s}{4}\left( 1-\tanh^2\left(\frac{k_s}{2}(u_{\tau}-\eta)\right) \right) (1-\dot{\tau})\right) \left(\operatorname{sgn}\left(\frac{\operatorname{sgn}(e_1)+1}{2}\right)\right) ke_{a_{\tau}} \\
        & +\alpha\left( \dot{u}_{\hat{\tau}} - \dot{u}_{\tau} \right) - \alpha\left(\operatorname{sgn}\left(\frac{\operatorname{sgn}(e_1)+1}{2}\right)\right) k e_a.
    \end{split}
\end{equation}
\normalsize
To better understand how $e_a$ propagates, we segregate terms in \eqref{rdot2} that can be upper bounded by a state-dependent function, and terms that can be upper bounded by a constant, such that

\footnotesize
\begin{equation} \label{rdot3}
    \begin{split}
        \dot{e}_a&= N_1 + N_2 - e_n -\alpha\left(\operatorname{sgn}\left(\frac{\operatorname{sgn}(e_1)+1}{2}\right)\right) k e_a+\alpha\left( \dot{u}_{\hat{\tau}} - \dot{u}_{\tau} \right) \\
        &+\left(\alpha-\frac{\beta k_s}{4}\left( 1-\tanh^2\left(\frac{k_s}{2}(u_{\tau}-\eta)\right) \right) (1-\dot{\tau})\right) \left(\operatorname{sgn}\left(\frac{\operatorname{sgn}(e_1)+1}{2}\right)\right) k e_{a_{\tau}}.
    \end{split}
\end{equation}
\normalsize
The boundable functions $N_1(t)\in\mathbb{R}$, $N_2(t)\in\mathbb{R}$ are defined as 
\small
\begin{align*}
    N_1&\triangleq -\dot{f}(\mathbf{x}_r)+x_r^{(n+1)},\\
    N_2&\triangleq -\dot{f}(\mathbf{x})+\dot{f}(\mathbf{x}_r)-\sum_{i=1}^{n-1}h_i^{(n-i+1)}(\mathbf{x})+\sum_{j=0}^{n-2}a_{n,j}e_1^{(j+2)}+\lambda\dot{e}_n+e_n ,
\end{align*}
\normalsize
where $\mathbf{x}_r(t)\triangleq \Big[ x_r,\dot{x}_r,...,\left(x_r^{(n-1)}\right)\Big]^T \in \mathbb{R}^n$.

\begin{remark} \label{remark1}
From Assumptions \ref{assumption_FunctionBound} and \ref{referenceassumption}, $N_1$ is upper bounded:
\begin{equation} \label{N1 bound}
    \sup_{t>0}\|N_1\|\leq c_1,
\end{equation}
where $c_1\in\mathbb{R}_+$ is a known positive constant.  
\end{remark} 
\begin{remark} \label{remark2}
Using Assumption \ref{assumption_FunctionBound} and Lemma 5 in \cite{kamalapurkar2013supporting}, $N_2$ is upper bounded as:
\begin{equation} \label{N2 bound}
    \|N_2\|\leq \|z\|\rho(\|z\|),
\end{equation}
where $\rho$ is a positive, radially unbounded, and strictly increasing function of $\|z\|$, and $z\in\mathbb{R}^{n+2}$ is a vector of error signals, $z\triangleq [e_1, e_2, \ldots,e_n,e_u,e_a]^T$.
\end{remark}

The Lyapunov-Krasovskii (LK) method extends the Lyapunov method to analyze the stability of differential equations with time delay~\cite{seuret2016d1}. This method selects energy functionals (functions of the system state) that are positive definite and decreasing, i.e., the derivative of the function is negative definite along the system trajectories. LK functionals are typically defined as sums of quadratic terms that depend on the delayed states~\cite{hetel2008equivalence}. We use LK-based functionals \cite{kolmanovskii2013introduction} similar to \cite{obuz2017unknown} for stability analysis. Let $ Q_1\in\mathbb{R}_+$, $Q_2\in\mathbb{R}_+$, $Q_3\in\mathbb{R}_+$ be
\begin{align}
    Q_1&= \omega_1\int_{t-\hat{\tau}}^t \|e_a(\theta)\|^2d\theta,\label{Q1} \\
    Q_2 &= \omega_2\int_{t-\tau}^t \|e_a(\theta)\|^2 d\theta, \label{Q2} \\
    Q_3 &= \omega_3\int_{t-(\Bar{\Tilde{\tau}}+\hat{\tau})}^t \int_{s}^t \|\dot{u}(\theta)\|^2d\theta ds, \label{Q3}
\end{align}
where $\omega_1$, $\omega_2$, and $\omega_3$ are positive constants. We define $\omega_2 \triangleq \frac{k\alpha}{1-\varphi_2}$. Additionally, define $y\in\mathbb{R}^{n+5}$ as
\begin{equation*} 
    y\triangleq [z, \sqrt{Q_1}, \sqrt{Q_2}, \sqrt{Q_3}]^T. 
\end{equation*}

Let $\varepsilon_1$ and $\varepsilon_2$ be positive constants. We will need auxiliary bounding constants $\sigma\in\mathbb{R}_+$, $\Delta\in\mathbb{R}_+$, which are

\footnotesize
\begin{align} \label{sigma bound}
    \sigma\triangleq &\min\bigg\{ 1,\left(1-\frac{\varepsilon_2}{2}\right),\left(\lambda-\left(\frac{\alpha}{2\varepsilon_1}+\frac{1}{2\varepsilon_2}\right) \right),\\
    &~~~~~~~~~~~~\left(\frac{\omega_3}{4\hat{\tau}}-\left(\frac{\alpha\varepsilon_1}{2}+\frac{k^2}{4\omega_1}\right)\right),\frac{k\alpha}{8} \bigg\}, \nonumber \\
    \Delta\triangleq &\frac{1}{2}\min\bigg\{\frac{\sigma}{2},\frac{\omega_3k^2}{4\omega_1}, \frac{\omega_3k^2}{4\omega_2},\frac{1}{4(\Bar{\Tilde{\tau}}+\hat{\tau})} \bigg\}. \label{delta bound}
\end{align}
\normalsize

We study stability of the error system over the domain  $\mathscr{D}_1\triangleq\Big\{y ~|~ \|y\|<\chi \Big\}$, where $\chi\triangleq\inf\Big\{\rho^{-1}(\xi)\Big\}$, $\forall\xi\in[\sqrt{\frac{\sigma k\alpha}{2}},\infty)$. 
Let $\Omega_\psi=\{\mu\in\mathbb{R}^n|z(\mu)<\psi\}$, then $\mathscr{D} \triangleq \mathscr{D}_1 \cap \Omega_\psi$ is the domain of attraction,  
i.e., every trajectory starting in $\mathscr{D}$ remains in $\mathscr{D}$ and approaches the origin as $t\rightarrow \infty$.


\begin{theorem} \label{theorem 1}
Under Assumptions \ref{assumption delay infinite state}--\ref{assumption error bound}, for the dynamics in \eqref{dymodel} with the controller in \eqref{controllereq} and \eqref{controller nu}, suppose that we select controller gains such that all of the following are satisfied:
\begin{equation} \label{gain conditions}
    \begin{split}
        \varepsilon_2 < 2, ~~~~
        \lambda >& \left(\frac{\alpha}{2\varepsilon_1}+\frac{1}{2\varepsilon_2}\right), ~~~~
        \omega_3 > 4\hat{\tau}\left(\frac{\alpha\varepsilon_1}{2}+\frac{k^2}{4\omega_1}\right), \\
        \Bar{\Tilde{\tau}} <& \left(\frac{\frac{k\alpha}{8}-\omega_1-\omega_2-\frac{\alpha}{2}-\omega_3 k \hat{\tau}}{\omega_3 k} \right).
    \end{split}
\end{equation}
Then, for the domain of attraction $\mathscr{D}$, the resulting closed-loop system is uniformly ultimately bounded,
\begin{equation}
    \limsup_{t\rightarrow\infty}\|e_1\|\leq \left( \sqrt{\frac{2c_1^2+k\alpha^2\Bar{\Tilde{\tau}}^2m^2}{k\alpha\Delta}} \right).\label{eq:UUBepsilon} 
\end{equation}
\end{theorem}
 
\noindent \textit{Proof.} See Appendix.

The tracking bound on the right side of \eqref{eq:UUBepsilon} is the $\epsilon$ of the reference tracking control problem that we set out to solve.

\section{Controller Application}\label{section5}
\subsection{Simulation Results}
To illustrate controller performance, we use the coagulation model \eqref{coagulationmodel} from  \cite{menezes2017targeted}. A state-space representation of the model that includes an unknown state-dependent power law input delay and input saturation as shown in \eqref{dymodel} is
\begin{equation} \label{sim coag model}
    \begin{split}
        &\dot{x}_1(t)=x_2(t)-d_1 x_1(t), \\
        &\dot{x}_2(t)=x_3(t) - d_2 x_2(t), \\
        &\dot{x}_3(t)=-d_3x_3(t)+g\left(u\left(t-\tau\left(x_3\left(t\right)\right)\right)\right).
    \end{split}
\end{equation}
Parameters of \eqref{sim coag model} for a trauma patient plasma sample are 
$d_1=1.1311$, $d_2=1.1362$, and $d_3=0.2727$, 
which we obtained from fits to experimental data. 
Reasonable saturation function parameters are: maximum value $\beta=50$, horizontal shift $\eta=75$, and growth rate $k_s=0.0224$. 

Time delay parameters observed from Fig.~\ref{fig1}A are $\gamma=4.48$, $k_d=0.322$. 
With these values, it is possible for the time delay magnitude to be the same order as the system time-constant. But to emphasize that the time delay can be unknown, we choose coarse estimates of the parameters underlying $\hat{\tau}$, the estimated time delay, as $\hat{\gamma}=1$ and $\hat{k_d}=1$. For illustration, we consider the initial condition $\mathbf{x_0}=[500,50,5]^T$. We choose controller gains based on \eqref{gain conditions} as $\lambda=0.1,\,\alpha=5,\,k=0.15$. 

We next present 
two clinically-relevant cases.

\begin{figure}[h!]
    \centering
    \includegraphics[width=.483\textwidth]{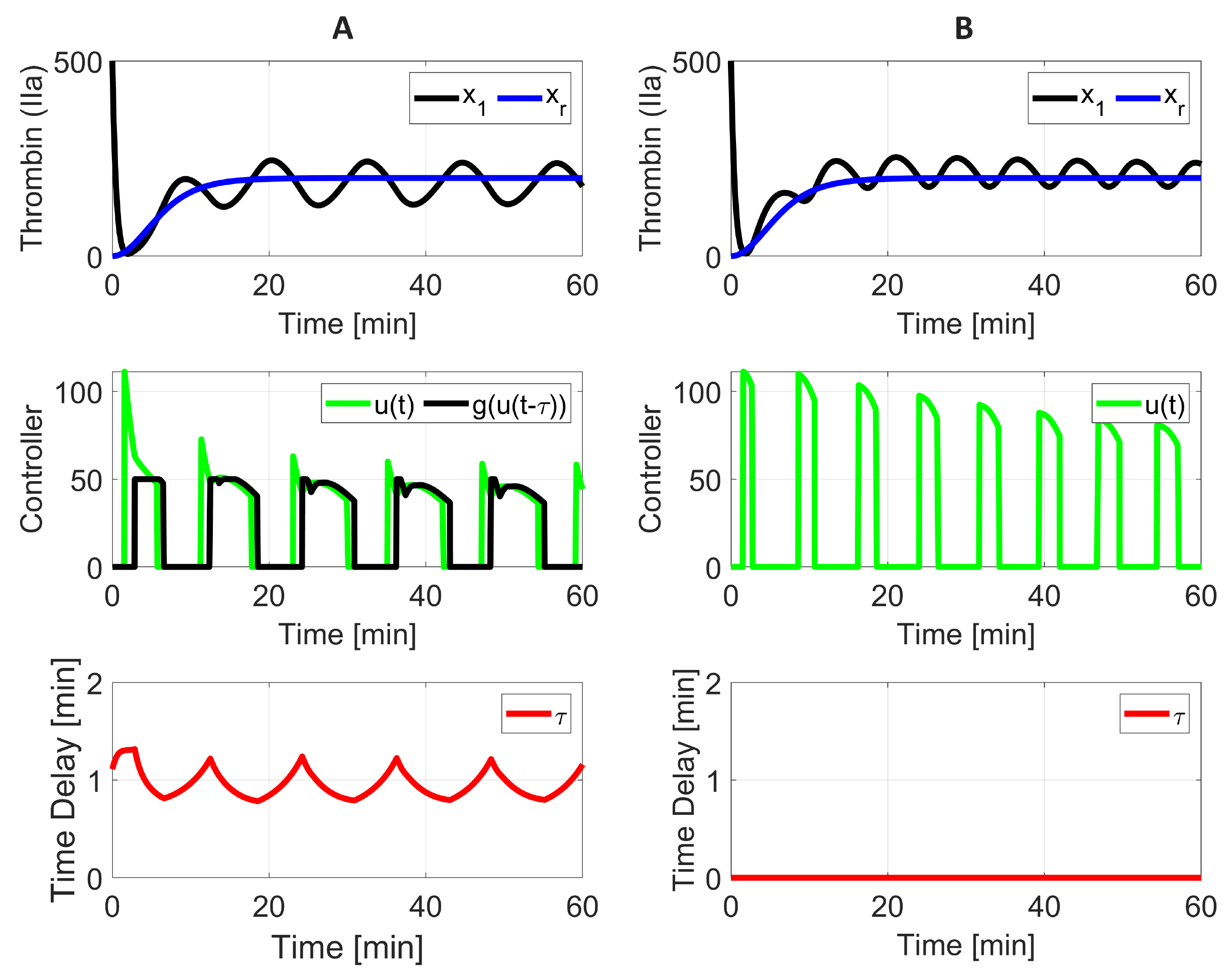}
    \vspace{-2em}
    \caption{Test of controller performance. (A) Satisfactory thrombin reference tracking (top panel) for Case 1.1, an initially-high thrombin level, with coagulation system input (middle panel), and associated time-varying input delay (bottom panel). (B) Our developed controller is also applicable to positive systems without input delay or saturation.}
    \label{fig_case1:1-2}
\end{figure}

\textbf{Case 1: Reduction of an initially-elevated thrombin level to a reference normal level, and maintenance at that level.} Many trauma patients experience high thrombin concentrations \cite{menezes2017targeted}. 
Unregulated concentrations of thrombin are the source of hypo- or hyper-coagulopathy, which can lead to bleeding, multi-organ failure, stroke, and death 
\cite{moore2021trauma}. 
Hence, there is a desire to regulate and maintain thrombin levels within a normal range. This case thus represents a recovery of thrombin concentration in injured humans. The reference that we wish to track is
\begin{equation*} \label{ref signal 1}
    x_{r_{1}}(t) = 200 \left(\tanh^2(0.15t)\right).
\end{equation*}
We present this case in the following four sub-cases:

\textbf{Case 1.1}. We investigate controller performance in the presence of the defined saturation and time delay. Fig.~\ref{fig_case1:1-2}A shows the satisfactory tracking performance for this case, the controller's periodic on-off inputs, and the associated delays. 

\textbf{Case 1.2}: We repeat Case 1.1 without the input saturation and the input time delay, Fig.~\ref{fig_case1:1-2}B. This case demonstrates the proposed controller's application for the class of positive systems that are not limited by the input nonlinearities studied in this paper. The depicted controller performance shows that there is only a small performance sacrifice that is made on reference tracking in traditional systems for an added benefit of reference tracking in more complex nonlinear systems.
\vspace{-1em}



\begin{figure}[h!]
    \centering
    \includegraphics[width=.483\textwidth]{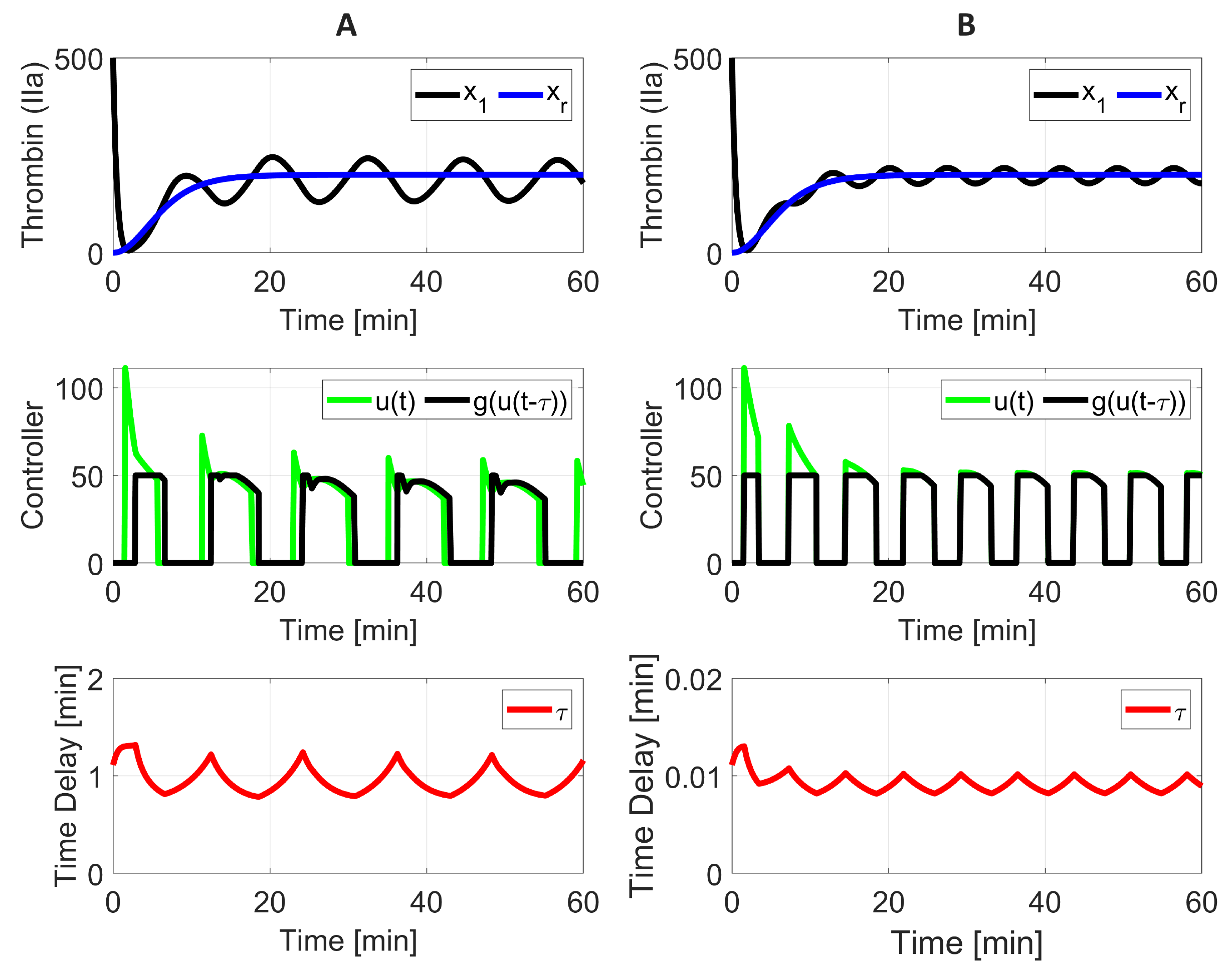}
    \vspace{-2em}
    \caption{The effects of time delay on controller performance. (A) The controller is unaffected by a poor time delay estimate, $\hat{\tau}$. (B) Satisfactory thrombin reference tracking (top panel) for Case 1.4, where the time delay magnitude is smaller than the system order (bottom panel). The coagulation system sees frequent inputs (middle panel).}
    \label{fig_case1:5-6}
    \vspace{-1em}
\end{figure}

\textbf{Case 1.3}: We show the lack of time delay estimation effect on controller performance. We reduce the time delay estimate magnitude $\hat{\tau}$ by one order, i.e., $\hat{\gamma}=0.1$, which is consequently a very poor estimate. Nevertheless, Fig.~\ref{fig_case1:5-6}A confirms satisfactory tracking and is identical to the results from Case 1.1 in Fig.~\ref{fig_case1:1-2}A.

\textbf{Case 1.4}: We show the time delay magnitude effects on controller performance. We reduced the time delay by two orders of magnitude, so $\gamma=0.0448$. The reference tracking in Fig.~\ref{fig_case1:5-6}B is improved compared to Case 1.1 due to the smaller time delay magnitude. 
This result confirms that, for applications where the time delay magnitude is not the same order as the system time-constant, the proposed controller still achieves tracking. Smaller time delays increase controller compensation frequency.

\begin{figure}[h!]
    \centering
    \includegraphics[width=.48\textwidth]{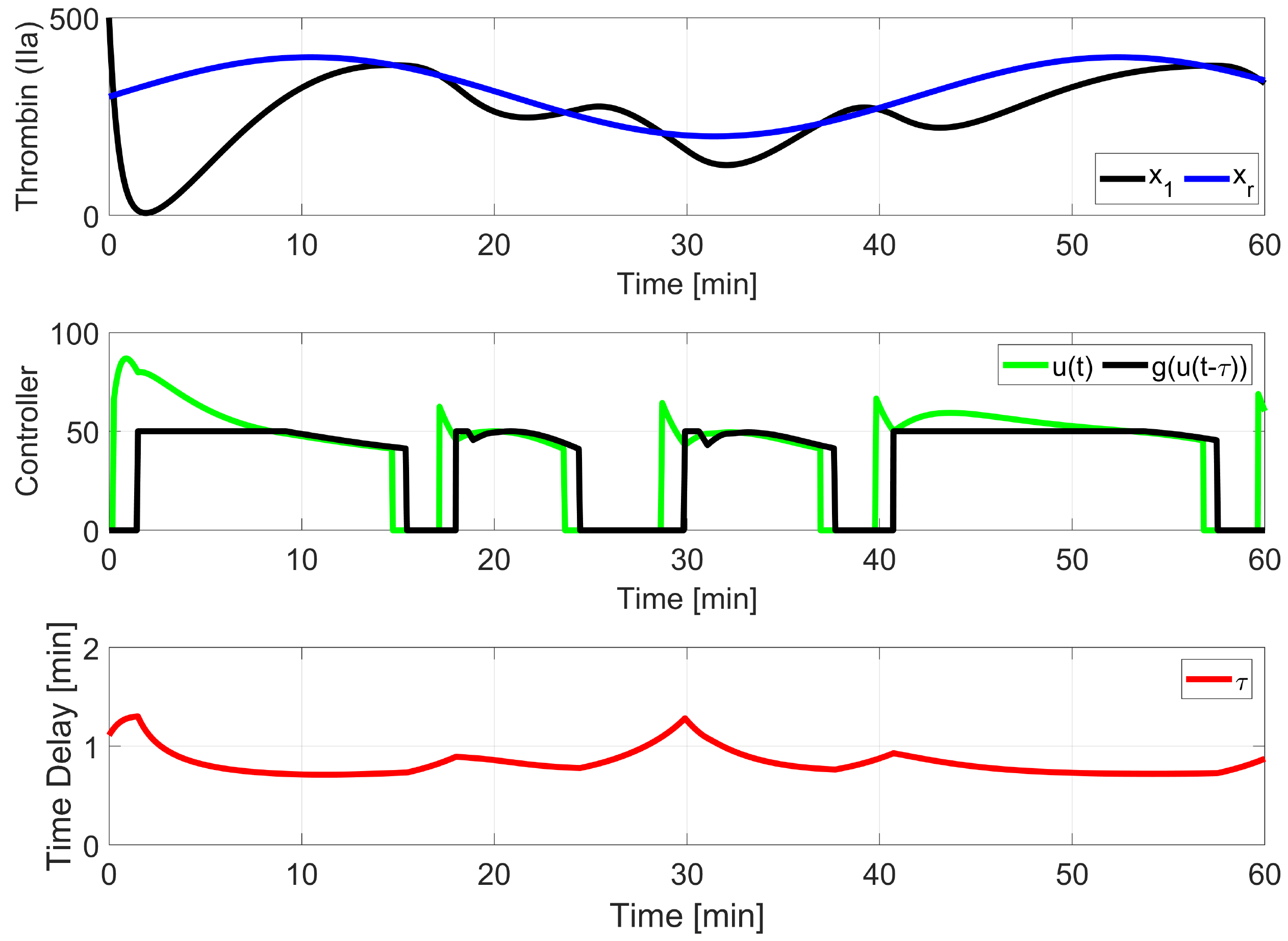}
    \vspace{-2em}
    \caption{Case 2: A time-varying sinusoidal reference to benefit patients 
    by leveraging both procoagulant and anticoagulant effects of thrombin. Satisfactory thrombin reference tracking (top panel), with coagulation system input (middle panel), and associated time-varying input delay (bottom panel).}
    \label{fig_case3}
\end{figure}

\textbf{Case 2: Tracking of a time varying reference.} We show controller robustness to variable references chosen based on patient condition. 
Since thrombin acts as both anticoagulant and procoagulant, a varying concentration can both reduce clotting and regulate bleeding. A reference signal for this case is \[x_{r_{2}}(t) = 100\Big(\sin(0.15t)\Big)+300.\]
Fig.~\ref{fig_case3} shows that the controller leverages pulse-width modulation of a saturated control signal to track $x_{r_{2}}$.
\vspace{-0.5em}

\vspace{-1.25em}
\subsection{Controller Performance Comparison}
\begin{figure}[h!]
    \centering
    \includegraphics[width=.483\textwidth]{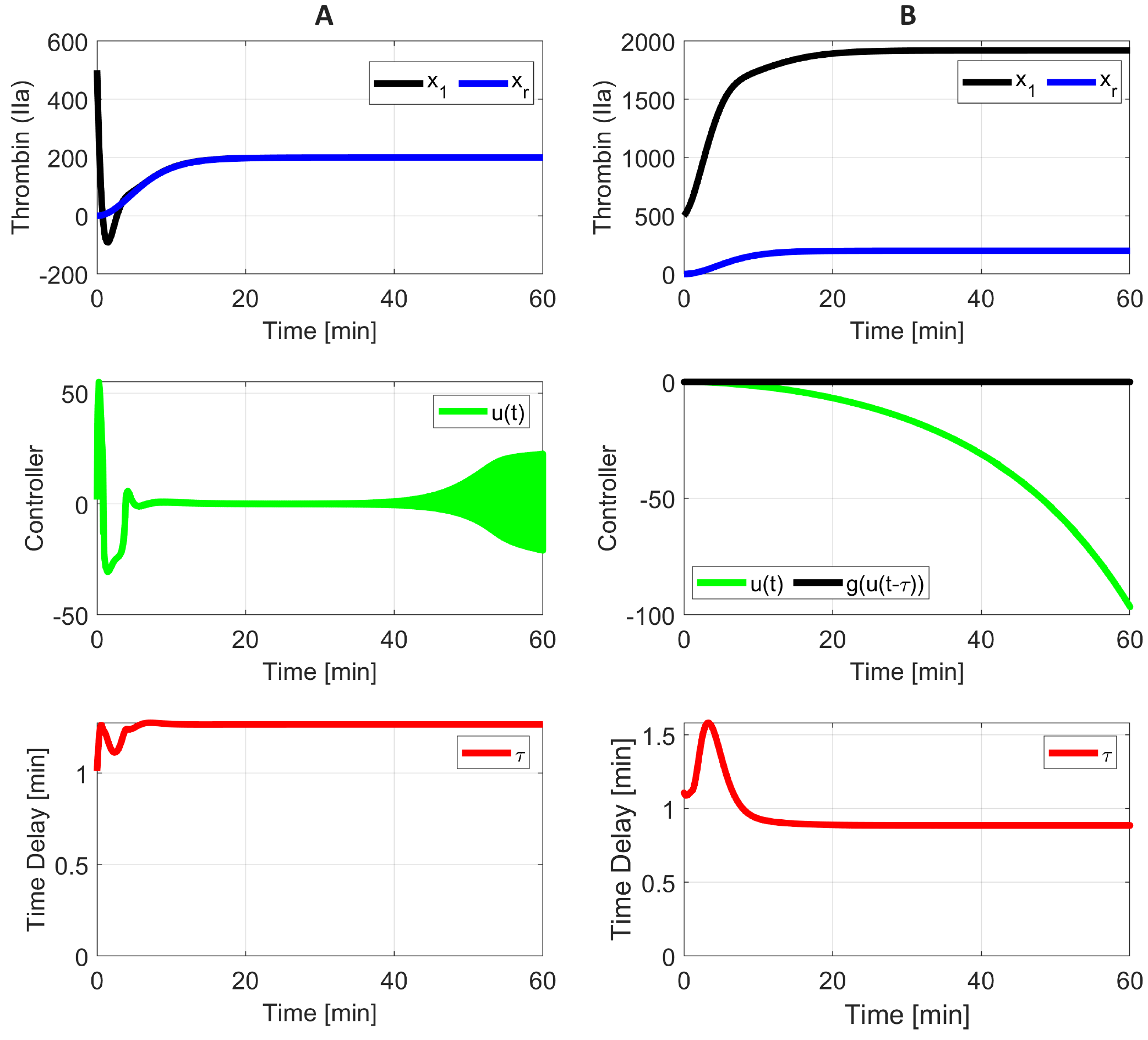}
    \vspace{-2em}
    \caption{Current appropriate state-of-the-art controllers in literature \cite{zhang2015exact} are not suitable for the proposed class of nonlinear systems in this study. (A) Satisfactory tracking is achieved confirming proper controller implementation. However, state $x_1$ (top panel) and control signal (middle panel) go negative, therefore violating the positive system definition. (B) Adding positive value enforcement to the same control law as A constrains tracking performance, and the control signal again goes negative. This figure confirms the necessity of this work.}
    \label{fig_controllerComparison}
\end{figure}

To contrast our controller performance against relevant state-of-the-art controllers, 
we repeated Case~1 using the controller in \cite{zhang2015exact}, which was developed for nonlinear systems with time delay and dead-zone input saturation. First, Fig.~\ref{fig_controllerComparison}A shows thrombin tracking, attesting to the proper simulation of the developed controller. However, both input $u(t)$ and thrombin level $x_1(t)$ go negative, violating the positive input and positive state requirements of positive systems. Moreover, when input saturation and positive state restrictions are applied to the simulation, thrombin tracking is not achieved, Fig.~\ref{fig_controllerComparison}B, and the controller signal violates the positive input requirement. Hence, the novel controller developed in this paper is necessary and beneficial. 

\vspace{-0.5em}

\vspace{-1.25em}
\section{Conclusions}\label{section6}
We have developed a satisfactorily-performing controller that can handle a class of systems restricted to positive states and inputs, cope with an uncertain state-dependent input delay that does not have to be of the same order as the system time-constant, and compensate for a control input that can saturate. 
Our simulation results confirm desired controller performance 
when applied to a pressing biomedical problem, 
which is treating trauma patients to remedy abnormal coagulation by personalizing their treatment through closed-loop protein delivery. 
Simulation results from different scenarios confirm that our results are also  applicable to positive systems that do not experience one or both of the input nonlinearities (input saturation and time delay) or do not have a good estimate of the time delay.  


\vspace{-0.3em}
\bibliography{References}

\begin{thebibliography}{10}
\expandafter\ifx\csname url\endcsname\relax
  \def\url#1{\texttt{#1}}\fi
\expandafter\ifx\csname urlprefix\endcsname\relax\def\urlprefix{URL }\fi
\expandafter\ifx\csname href\endcsname\relax
  \def\href#1#2{#2} \def\path#1{#1}\fi

\bibitem{farina2011positive}
L.~Farina, S.~Rinaldi, Positive Linear Systems: Theory and Applications,
  Vol.~50, John Wiley \& Sons, 2011.

\bibitem{kaczorek2012positive}
T.~Kaczorek, Positive 1D and 2D Systems, Springer Science \& Business Media,
  2012.

\bibitem{blanchini2015switched}
F.~Blanchini, P.~Colaneri, M.~E. Valcher, Switched positive linear systems,
  Foundations and Trends in Systems and Control 2~(2) (2015) 101--273.

\bibitem{rantzer2018tutorial}
A.~Rantzer, M.~E. Valcher, A tutorial on positive systems and large scale
  control, in: 2018 IEEE Conference on Decision and Control (CDC), IEEE, 2018,
  pp. 3686--3697.

\bibitem{rantzer2021scalable}
A.~Rantzer, M.~E. Valcher, Scalable control of positive systems, Annual Review
  of Control, Robotics, and Autonomous Systems 4 (2021) 319--341.

\bibitem{shorten2006positive}
R.~Shorten, F.~Wirth, D.~Leith, A positive systems model of {TCP}-like
  congestion control: Asymptotic results, IEEE/ACM Transactions on Networking
  14~(3) (2006) 616--629.

\bibitem{murray2002mathematical}
J.~D. Murray, Mathematical Biology, Volume {I}, An Introduction, Springer
  Verlag, New York, 2002.

\bibitem{haddad2010nonnegative}
W.~M. Haddad, V.~Chellaboina, Q.~Hui, Nonnegative and Compartmental Dynamical
  Systems, Princeton University Press, 2010.

\bibitem{blanchini2014piecewise}
F.~Blanchini, G.~Giordano, Piecewise-linear {L}yapunov functions for structural
  stability of biochemical networks, Automatica 50~(10).

\bibitem{briat2016antithetic}
C.~Briat, A.~Gupta, M.~Khammash, Antithetic integral feedback ensures robust
  perfect adaptation in noisy biomolecular networks, Cell Systems 2~(1) (2016)
  15--26.

\bibitem{hovorka2004nonlinear}
R.~Hovorka, V.~Canonico, L.~J. Chassin, U.~Haueter, M.~Massi-Benedetti,
  M.~O.~F. Federici, T.~R. Pieber, H.~C. Schaller, L.~Schaupp, T.~Vering, M.~E.
  Wilinska, Nonlinear model predictive control of glucose concentration in
  subjects with type 1 diabetes, Physiological Measurement 25~(4).

\bibitem{shu2008positive}
Z.~Shu, J.~Lam, H.~Gao, B.~Du, L.~Wu, Positive observers and dynamic
  output-feedback controllers for interval positive linear systems, IEEE
  Transactions on Circuits and Systems I: Regular Papers 55~(10).

\bibitem{shen2016static}
J.~Shen, J.~Lam, Static output-feedback stabilization with optimal {L}1-gain
  for positive linear systems, Automatica 63 (2016) 248--253.

\bibitem{coxson1987positive}
P.~G. Coxson, H.~Shapiro, Positive input reachability and controllability of
  positive systems, Linear Algebra and its Applications 94 (1987) 35--53.

\bibitem{benvenuti2004tutorial}
L.~Benvenuti, L.~Farina, A tutorial on the positive realization problem, IEEE
  Transactions on Automatic Control 49~(5) (2004) 651--664.

\bibitem{valcher2009reachability}
M.~E. Valcher, Reachability properties of continuous-time positive systems,
  IEEE Transactions on Automatic Control 54~(7) (2009) 1586--1590.

\bibitem{guiver2014positive}
C.~Guiver, D.~Hodgson, S.~Townley, Positive state controllability of positive
  linear systems, Systems \& Control Letters 65 (2014) 23--29.

\bibitem{eden2016positive}
J.~Eden, Y.~Tan, D.~Lau, D.~Oetomo, On the positive output controllability of
  linear time invariant systems, Automatica 71 (2016) 202--209.

\bibitem{briat2020biology}
C.~Briat, A biology-inspired approach to the positive integral control of
  positive systems: The antithetic, exponential, and logistic integral
  controllers, SIAM Journal on Applied Dynamical Systems 19~(1).

\bibitem{wu2010stability}
M.~Wu, Y.~He, J.-H. She, Stability Analysis and Robust Control of Time-Delay
  Systems, Vol.~22, Springer, 2010.

\bibitem{otto2019nonlinear}
A.~Otto, W.~Just, G.~Radons, Nonlinear dynamics of delay systems: An overview,
  Philosophical Transactions of Royal Society A 377.

\bibitem{hale1993introduction}
J.~K. Hale, S.~M.~V. Lunel, L.~S. Verduyn, S.~M.~V. Lunel, Introduction to
  Functional Differential Equations, Vol.~99, Springer Science \& Business
  Media, 1993.

\bibitem{park2014stability}
P.~Park, W.~I. Lee, S.~Y. Lee, Stability on time delay systems: A survey,
  Journal of Institute of Control, Robotics and Systems 20~(3).

\bibitem{li1997delay}
X.~Li, C.~E. De~Souza, Delay-dependent robust stability and stabilization of
  uncertain linear delay systems: A linear matrix inequality approach, IEEE
  Transactions on Automatic Control 42~(8) (1997) 1144--1148.

\bibitem{sheng2021switched}
Z.~Sheng, Z.~Sun, V.~Molazadeh, N.~Sharma, Switched control of an
  n-degree-of-freedom input delayed wearable robotic system, Automatica 125
  (2021) 109455.

\bibitem{obuz2017unknown}
S.~Obuz, J.~R. Klotz, R.~Kamalapurkar, W.~Dixon, Unknown time-varying input
  delay compensation for uncertain nonlinear systems, Automatica 76 (2017)
  222--229.

\bibitem{bekiaris2012compensation}
N.~Bekiaris-Liberis, M.~Krstic, Compensation of state-dependent input delay for
  nonlinear systems, IEEE Transactions on Automatic Control 58~(2) (2012)
  275--289.

\bibitem{alibeji2017pid}
N.~Alibeji, N.~Sharma, A pid-type robust input delay compensation method for
  uncertain euler--lagrange systems, IEEE Transactions on Control Systems
  Technology 25~(6) (2017) 2235--2242.

\bibitem{deng2021state}
Y.~Deng, V.~L{\'e}chapp{\'e}, E.~Moulay, F.~Plestan, State feedback control and
  delay estimation for {LTI} system with unknown input-delay, International
  Journal of Control 94~(9) (2021) 2369--2378.

\bibitem{koo2020output}
M.-S. Koo, H.-L. Choi, Output feedback regulation of a class of high-order
  feedforward nonlinear systems with unknown time-varying delay in the input
  under measurement sensitivity, International Journal of Robust and Nonlinear
  Control 30~(12) (2020) 4744--4763.

\bibitem{nguyen2021state}
C.~M. Nguyen, C.~P. Tan, H.~Trinh, State and delay reconstruction for nonlinear
  systems with input delays, Applied Mathematics and Computation 390 (2021)
  125609.

\bibitem{cai2018adaptive}
J.~Cai, J.~Wan, H.~Que, Q.~Zhou, L.~Shen, Adaptive actuator failure
  compensation control of second-order nonlinear systems with unknown time
  delay, IEEE Access 6 (2018) 15170--15177.

\bibitem{zhang2015exact}
Z.~Zhang, S.~Xu, B.~Zhang, Exact tracking control of nonlinear systems with
  time delays and dead-zone input, Automatica 52 (2015) 272--276.

\bibitem{perez2009saturation}
N.~O. P{\'e}rez-Arancibia, T.-C. Tsao, J.~S. Gibson, Saturation-induced
  instability and its avoidance in adaptive control of hard disk drives, IEEE
  Transactions on Control Systems Technology 18~(2) (2009) 368--382.

\bibitem{rami2007positive}
M.~A. Rami, U.~Helmke, F.~Tadeo, Positive observation problem for linear
  time-delay positive systems, in: 2007 Mediterranean Conference on Control \&
  Automation, IEEE, 2007, pp. 1--6.

\bibitem{liu2010stability}
X.~Liu, W.~Yu, L.~Wang, Stability analysis for continuous-time positive systems
  with time-varying delays, IEEE Transactions on Automatic Control 55~(4)
  (2010) 1024--1028.

\bibitem{liu2009constrained}
X.~Liu, Constrained control of positive systems with delays, IEEE Transactions
  on Automatic Control 54~(7) (2009) 1596--1600.

\bibitem{zhu2013exponential}
S.~Zhu, M.~Meng, C.~Zhang, Exponential stability for positive systems with
  bounded time-varying delays and static output feedback stabilization, Journal
  of the Franklin Institute 350~(3) (2013) 617--636.

\bibitem{zhang2015h}
Q.~Zhang, Y.~Zhang, B.~Du, T.~Tanaka, H$\infty$ control via dynamic output
  feedback for positive systems with multiple delays, IET Control Theory \&
  Applications 9~(17) (2015) 2574--2580.

\bibitem{hong2019optimization}
L.~V. Hien, M.~T. Hong, An optimization approach to static output-feedback
  control of {LTI} positive systems with delayed measurements, Journal of the
  Franklin Institute 356~(10) (2019) 5087--5103.

\bibitem{arogbonlo2019design}
V.~T. Huynh, A.~Arogbonlo, H.~Trinh, A.~M.~T. Oo, Design of observers for
  positive systems with delayed input and output information, IEEE Transactions
  on Circuits and Systems II: Express Briefs 67~(1) (2020) 107--111.

\bibitem{ajmeri2015simple}
M.~Ajmeri, A.~Ali, Simple tuning rules for integrating processes with large
  time delay, Asian Journal of Control 17~(5) (2015) 2033--2040.

\bibitem{klamka1996constrained}
J.~Klamka, Constrained controllability of nonlinear systems, Journal of
  Mathematical Analysis and Applications 201~(2) (1996) 365--374.

\bibitem{naim2018controllability}
M.~Naim, F.~Lahmidi, A.~Namir, Controllability and observability analysis of
  nonlinear positive discrete systems, Discrete Dynamics in Nature and Society
  2018.

\bibitem{greenstreet1999reachability}
M.~R. Greenstreet, I.~Mitchell, Reachability analysis using polygonal
  projections, in: International Workshop on Hybrid Systems: Computation and
  Control, Springer, 1999, pp. 103--116.

\bibitem{asarin2003reachability}
E.~Asarin, T.~Dang, A.~Girard, Reachability analysis of nonlinear systems using
  conservative approximation, in: International Workshop on Hybrid Systems:
  Computation and Control, 2003, pp. 20--35.

\bibitem{liu2015stability}
X.~Liu, Stability analysis of a class of nonlinear positive switched systems
  with delays, Nonlinear Analysis: Hybrid Systems 16 (2015) 1--12.

\bibitem{sabatier2020power}
J.~Sabatier, Power law type long memory behaviors modeled with distributed time
  delay systems, Fractal and Fractional 4~(1) (2020) 1.

\bibitem{sun2004explanation}
K.~Sun, Explanation of log-normal distributions and power-law distributions in
  biology and social science, Tech. rep., University of Illinois at
  Urbana-Champaign, Department of Physics (2004).

\bibitem{meyer1994bi}
P.~Meyer, Bi-logistic growth, Technological Forecasting and Social Change
  47~(1) (1994) 89--102.

\bibitem{csimcsek2019power}
E.~{\c{S}}im{\c{s}}ek, M.~Kim, Power-law tail in lag time distribution
  underlies bacterial persistence, Proceedings of the National Academy of
  Sciences 116~(36) (2019) 17635--17640.

\bibitem{menezes2017targeted}
A.~A. Menezes, R.~F. Vilardi, A.~P. Arkin, M.~J. Cohen, Targeted clinical
  control of trauma patient coagulation through a thrombin dynamics model,
  Science Translational Medicine 9~(371).

\bibitem{narayanan1999multifunctional}
S.~Narayanan, Multifunctional roles of thrombin, Annals of Clinical \&
  Laboratory Science 29~(4) (1999) 275--280.

\bibitem{luenberger1979introduction}
D.~G. Luenberger, Introduction to Dynamic Systems: Theory Models and
  Applications, Wiley, 1979.

\bibitem{ghetmiri2021personalized}
D.~E. Ghetmiri, M.~J. Cohen, A.~A. Menezes, Personalized modulation of
  coagulation factors using a thrombin dynamics model to treat trauma-induced
  coagulopathy, npj Systems Biology and Applications 7 (2021) 44.

\bibitem{gonzalez2016trauma}
E.~Gonzalez, H.~B. Moore, E.~E. Moore, Trauma Induced Coagulopathy, Springer,
  2016.

\bibitem{ghetmiri2021controltheoretic}
D.~E. Ghetmiri, M.~E. Perez~Blanco, M.~J. Cohen, A.~A. Menezes,
  Control-theoretic modeling and prediction of blood clot viscoelasticity in
  trauma patients, Proceedings of the Modeling, Estimation and Control
  Conference (MECC 2021), IFAC-PapersOnLine 54~(20) (2021) 232--237.

\bibitem{ghetmiri2022nonlinear}
D.~E. Ghetmiri, A.~A. Menezes, Nonlinear dynamic modeling and model predictive
  control of thrombin generation to treat trauma-induced coagulopathy,
  International Journal of Robust and Nonlinear Control (2022) rnc.5963.

\bibitem{sleet2011injury}
D.~A. Sleet, L.~L. Dahlberg, S.~V. Basavaraju, J.~A. Mercy, L.~C. McGuire,
  A.~Greenspan, et~al., Injury prevention, violence prevention, and trauma
  care: Building the scientific base, MMWR Surveill Summ 60~(Suppl 4) (2011)
  78--85.

\bibitem{moore2021trauma}
E.~E. Moore, H.~B. Moore, L.~Z. Kornblith, M.~D. Neal, M.~Hoffman, N.~J. Mutch,
  H.~Sch{\"o}chl, B.~J. Hunt, A.~Sauaia, Trauma-induced coagulopathy, Nature
  Reviews Disease Primers 7~(1) (2021) 1--23.

\bibitem{cantle2017prediction}
P.~M. Cantle, B.~A. Cotton, Prediction of massive transfusion in trauma,
  Critical Care Clinics 33~(1) (2017) 71--84.

\bibitem{allaby2013dictionary}
M.~Allaby, A dictionary of geology and earth sciences, Oxford University Press,
  2013.

\bibitem{young2017dynamics}
J.~T. Young, T.~S. Hatakeyama, K.~Kaneko, Dynamics robustness of cascading
  systems, PLoS Computational Biology 13~(3) (2017) e1005434.

\bibitem{guo2017critical}
H.~Guo, C.~Zheng, H.~H.-C. Iu, T.~Fernando, A critical review of cascading
  failure analysis and modeling of power system, Renewable and Sustainable
  Energy Reviews 80 (2017) 9--22.

\bibitem{kaczorek2015analysis}
T.~Kaczorek, Analysis of positivity and stability of discrete-time and
  continuous-time nonlinear systems, Computational Problems of Electrical
  Engineering 5~(1) (2015) 11--16.

\bibitem{xian2004continuous}
B.~Xian, D.~M. Dawson, M.~S. de~Queiroz, J.~Chen, A continuous asymptotic
  tracking control strategy for uncertain nonlinear systems, IEEE Transactions
  on Automatic Control 49~(7) (2004) 1206--1211.

\bibitem{xian2004continuous_incollection}
B.~Xian, M.~S. De~Queiroz, D.~M. Dawson, A continuous control mechanism for
  uncertain nonlinear systems, in: Optimal Control, Stabilization and Nonsmooth
  Analysis, Springer, 2004, pp. 251--264.

\bibitem{sole1992structural}
R.~V. Sol{\'e}, J.~Valls, On structural stability and chaos in biological
  systems, Journal of Theoretical Biology 155~(1) (1992) 87--102.

\bibitem{bracewell1986fourier}
R.~N. Bracewell, R.~N. Bracewell, The Fourier transform and its applications,
  Vol. 31999, McGraw-Hill New York, 1986.

\bibitem{kamalapurkar2013supporting}
R.~Kamalapurkar, J.~A. Rosenfeld, J.~Klotz, R.~J. Downey, W.~E. Dixon,
  Supporting lemmas for rise-based control methods, arXiv preprint
  arXiv:1306.3432.

\bibitem{seuret2016d1}
A.~Seuret, F.~Gouaisbaut, L.~Baudouin, Overview of {L}yapunov methods for
  time-delay systems, Ph.D. thesis, Laboratory for Analysis and Architecture of
  Systems (2016).

\bibitem{hetel2008equivalence}
L.~Hetel, J.~Daafouz, C.~Iung, Equivalence between the {L}yapunov-{K}rasovskii
  functionals approach for discrete delay systems and that of the stability
  conditions for switched systems, Nonlinear Analysis: Hybrid Systems 2~(3)
  (2008) 697--705.

\bibitem{kolmanovskii2013introduction}
V.~Kolmanovskii, A.~Myshkis, Introduction to the theory and applications of
  functional differential equations, Vol. 463, Springer Science \& Business
  Media, 2013.

\bibitem{cortes2008discontinuous}
J.~Cortes, Discontinuous dynamical systems, IEEE Control Systems Magazine
  28~(3) (2008) 36--73.

\bibitem{filippov2013differential}
A.~F. Filippov, Differential equations with discontinuous righthand sides:
  control systems, Vol.~18, Springer Science \& Business Media, 2013.

\end{thebibliography}
\flushcolsend
\newpage
\section*{Appendix}
\noindent\textbf{Theorem 1.}

\noindent \textit{Proof.} We prove this theorem directly, using a Lyapunov stability analysis. 

Let $V$ be a Lyapunov function candidate defined as
\begin{equation} \label{lyapunov func}
    V \triangleq \frac{1}{2}\sum_{i=1}^n e_i^2 + \frac{1}{2}e_a^2 +\frac{1}{2}e_u^2 + Q_1+Q_2+Q_3,
\end{equation}
where $V>0,\,\forall t>0$. This definition is such that
\begin{equation} \label{lyap ineq}
    \frac{1}{2}\|y\|^2\leq V(y) \leq \|y\|^2. 
\end{equation}

The derivative of the first term in \eqref{lyapunov func} can be obtained using \eqref{error sys en}--\eqref{enerror}, \eqref{general ei}  with $i=n$, and \eqref{e_a error} as
\begin{equation} \label{sum error der}
    \begin{split}
        \sum_{i=1}^n e_i\dot{e}_i&=e_n\dot{e}_n+\sum_{i=1}^{n-1}e_i\dot{e}_i=e_n\dot{e}_n+e_{n-1}e_n-\sum_{i=1}^{n-1}e_i^2, \\
        &= e_n e_a-\lambda e_n^2-\alpha e_n e_u + e_{n-1}e_n-\sum_{i=1}^{n-1}e_i^2.
    \end{split}
\end{equation}
Using the Leibniz Rule, we can obtain the derivatives of \eqref{Q1}-\eqref{Q3} as
\begin{align} 
    \dot{Q}_1&=\omega_1(\|e_a\|^2-\|e_{a_{\hat{\tau}}}\|^2),\label{Q1der}\\
    \dot{Q}_2&=\omega_2(\|e_a\|^2-(1-\dot{\tau})\|e_{a_{\tau}}\|^2), \label{Q2der}\\
    \dot{Q}_3&=\omega_3\left( (\Bar{\Tilde{\tau}}+\hat{\tau})\|\dot{u}(t)\|^2 -  \int_{t-(\Bar{\Tilde{\tau}}+\hat{\tau})}^t\| \dot{u}(\theta) \|^2 d\theta \right). \label{Q3der}
\end{align}

With \eqref{eu error}, \eqref{e_a error}, \eqref{eq:doteu},  \eqref{udot}, \eqref{rdot3}, \eqref{sum error der}, and \eqref{Q1der}--\eqref{Q3der}, the derivative of \eqref{lyapunov func} is

\footnotesize
\begin{equation} \label{vdot 1}
    \begin{split}
        \dot{V}&= e_n e_a-\lambda e_n^2-\alpha e_n e_u + e_{n-1}e_n-\sum_{i=1}^{n-1}e_i^2 \\
        &+e_a \bigg[ -\alpha\left(\operatorname{sgn}\left(\frac{\operatorname{sgn}(e_1)+1}{2}\right)\right) k e_a+\alpha\left( \dot{u}_{\hat{\tau}} - \dot{u}_{\tau} \right) \\
        &+ \left(\alpha-\frac{\beta k_s}{4}\left( 1-\tanh^2\left(\frac{k_s}{2}(u_{\tau}-\eta)\right) \right) (1-\dot{\tau})\right) \left(\operatorname{sgn}\left(\frac{\operatorname{sgn}(e_1)+1}{2}\right)\right) k e_{a_{\tau}}\\
        &  +N_1 + N_2 - e_n \bigg] + e_u\left(\operatorname{sgn}\left(\frac{\operatorname{sgn}(e_1)+1}{2}\right)\right)\left(k e_{a_{\hat{\tau}}}-k e_a\right)\\
        & +\omega_1(\|e_a\|^2-\|e_{a_{\hat{\tau}}}\|^2) +\omega_2(\|e_a\|^2-(1-\dot{\tau})\|e_{a_{\tau}}\|^2)\\ 
        & +\omega_3\left( (\Bar{\Tilde{\tau}}+\hat{\tau})k^2\|e_a\|^2 -  \int_{t-(\Bar{\Tilde{\tau}}+\hat{\tau})}^t\| \dot{u}(\theta) \|^2 d\theta \right).
    \end{split}
\end{equation}
\normalsize

We will need the following small lemma, Lemma \ref{lem:delayassumptionbound}, to upper bound the $(1-\dot{\tau})$ term. The
purpose of Lemma \ref{lem:delayassumptionbound} is to show that there exists a $\varphi_2$ bounding $|\dot{\tau}|$, thus to show that Assumption~\ref{delayassumption} holds. As highlighted in the main text after Assumption~\ref{delayassumption}, $\varphi_2$ is an arbitrarily large number that can depend on the arbitrarily small $\varphi_1$, such that it is a non-restrictive bound on how fast the delay varies, $|\dot{\tau}|< \varphi_2$ $\forall t\geq 0$.

\begin{lemma}\label{lem:delayassumptionbound}
For any $\varphi_2>d_nk_d\gamma\varphi_1^{-k_d}$, we have $(1-\dot{\tau}) > (1 - \varphi_2)$.
\end{lemma}

\begin{proof}
We prove this lemma directly.  Consider that
\begin{align*}
    \dot{\tau} &= -k_d\gamma x_n^{-(k_d +1)} \dot{x}_n,\\
    &= -k_d\gamma x_n^{-(k_d +1)} (f(\mathbf{x})+g(u)),\\
    &= -k_d\gamma x_n^{-(k_d +1)} (-d_nx_n + \mathcal{F}(\mathbf{x})+g(u)),\\
    &= d_n k_d\gamma x_n^{-k_d} -k_d\gamma x_n^{-(k_d +1)} ( \mathcal{F}(\mathbf{x})+g(u)).
\end{align*}
Since $\mathcal{F}(\mathbf{x})$ and $g(u)$ are positive functions,
\begin{align*}
    \dot{\tau} &< d_n k_d\gamma x_n^{-k_d},\\
    &< d_n k_d\gamma \varphi_1^{-k_d},\\
    &< \varphi_2.
\end{align*}
Accordingly, $(1-\dot{\tau}) > (1 - \varphi_2)$. 
\end{proof}

After canceling common terms, and using 
Lemma \ref{lem:delayassumptionbound} 
and the fact that $\left(\operatorname{sgn}\left(\frac{\operatorname{sgn}(e_1)+1}{2} \right)\right)\leq 1$, \eqref{vdot 1} can be upper bounded as

\small
\begin{equation} \label{vdot 2}
    \begin{split}
        \dot{V} \leq &-\sum_{i=1}^{n-1}\|e_i\|^2 -\lambda e_n^2-\alpha |e_n e_u|+ |e_{n-1}e_n|\\
        &  + k\left(\alpha-\frac{\beta k_s}{4}\left( 1-\tanh^2\left(\frac{k_s}{2}(u_{\tau}-\eta)\right) \right) (1-\varphi_2)\right)|e_a e_{a_{\tau}}| \\
        & -\alpha k \|e_a\|^2 + \alpha|e_a\left( \dot{u}_{\hat{\tau}} - \dot{u}_{\tau} \right)| + e_a(N_1+N_2)\\
        & +k|e_u e_{a_{\hat{\tau}}}|-k|e_u e_a| \\
        & +\omega_1(\|e_a\|^2-\|e_{a_{\hat{\tau}}}\|^2) +\omega_2(\|e_a\|^2-(1-\varphi_2)\|e_{a_{\tau}}\|^2)\\ 
        & +\omega_3\left( (\Bar{\Tilde{\tau}}+\hat{\tau})k^2\|e_a\|^2 -  \int_{t-(\Bar{\Tilde{\tau}}+\hat{\tau})}^t\| \dot{u}(\theta) \|^2 d\theta \right).
    \end{split}
\end{equation}
\normalsize

Since
\begin{align*}
    0 \leq &1-\tanh^2(\cdot), \\
    \Rightarrow 0 \geq & -\frac{\beta k_s}{4}(1-\tanh^2(\cdot))(1-\varphi_2), \\
    \Rightarrow \alpha \geq & \alpha-\frac{\beta k_s}{4}(1-\tanh^2(\cdot))(1-\varphi_2),
\end{align*}
then \eqref{vdot 2} can be upper bounded as 
\begin{equation} \label{vdot 3}
    \begin{split}
        \dot{V} \leq &-\sum_{i=1}^{n-1}\|e_i\|^2 -\lambda e_n^2-\alpha |e_n e_u|+ |e_{n-1}e_n|  + k\alpha|e_a e_{a_{\tau}}| \\
        & -\alpha k \|e_a\|^2 + \alpha|e_a\left( \dot{u}_{\hat{\tau}} - \dot{u}_{\tau} \right)| + e_a(N_1+N_2)\\
        & +k|e_u e_{a_{\hat{\tau}}}|-k|e_u e_a| \\
        & +\omega_1(\|e_a\|^2-\|e_{a_{\hat{\tau}}}\|^2) +\omega_2(\|e_a\|^2-(1-\varphi_2)\|e_{a_{\tau}}\|^2)\\ 
        & +\omega_3\left( (\Bar{\Tilde{\tau}}+\hat{\tau})k^2\|e_a\|^2 -  \int_{t-(\Bar{\Tilde{\tau}}+\hat{\tau})}^t\| \dot{u}(\theta) \|^2 d\theta \right).
    \end{split}
\end{equation}

Using Young?s Inequality, the following inequalities can be obtained
\begin{align}\label{ineq 1}
    &|e_n e_u| \leq \frac{1}{2\varepsilon_1}\|e_n\|^2+\frac{\varepsilon_1}{2}\|e_u\|^2, \\ 
    &|e_{n-1}e_n| \leq \frac{\varepsilon_2}{2}\|e_{n-1}\|^2 +\frac{1}{2\varepsilon_2}\|e_n\|^2, \\ 
    &|e_a( \dot{u}_{\hat{\tau}} - \dot{u}_{\tau} )| \leq \frac{1}{2}\|e_a\|^2+ \frac{1}{2}\|\dot{u}_{\hat{\tau}} - \dot{u}_{\tau}\|^2. \label{ineq 3}
\end{align}
By using \eqref{ineq 1}--\eqref{ineq 3}, and Remarks \ref{remark1} and \ref{remark2}, \eqref{vdot 3} can be written as

\footnotesize
\begin{equation} \label{vdot 4}
    \begin{split}
        \dot{V} \leq &-\sum_{i=1}^{n-2}\|e_i\|^2 -\left(1-\frac{\varepsilon_2}{2}\right)\|e_{n-1}\|^2  -\left(\lambda-\left(\frac{\alpha}{2\varepsilon_1}+\frac{1}{2\varepsilon_2}\right) \right) \|e_n\|^2\\ &+\frac{\alpha\varepsilon_1}{2}\|e_u\|^2 + k\alpha\|e_a\|\|e_{a_{\tau}}\| +k\|e_u\| \|e_{a_{\hat{\tau}}}\| -k\|e_u e_a\|\\
        & -\alpha k \|e_a\|^2 + \frac{\alpha}{2}\|e_a\|^2+ \frac{\alpha}{2}\|\dot{u}_{\hat{\tau}} - \dot{u}_{\tau}\|^2 + c_1\|e_a\|+\|e_a\|\|z\|\rho(\|z\|)\\
        & +\omega_1\left(\|e_a\|^2-\|e_{a_{\hat{\tau}}}\|^2\right) +\omega_2\left(\|e_a\|^2-(1-\varphi_2)\|e_{a_{\tau}}\|^2\right)\\ 
        & +\omega_3\left( (\Bar{\Tilde{\tau}}+\hat{\tau})k^2\|e_a\|^2 -  \int_{t-(\Bar{\Tilde{\tau}}+\hat{\tau})}^t\| \dot{u}(\theta) \|^2 d\theta \right).
    \end{split}
\end{equation}
\normalsize

By completing the squares, we can develop two inequalities. First:

\footnotesize
\begin{equation} \label{ineq 4}
    \begin{split}
        &k\alpha\|e_a\|\|e_{a_{\tau}}\|-\omega_2(1-\varphi_2)\|e_{a_{\tau}}\|^2, \\
        = & -\omega_2(1-\varphi_2)\left(\|e_{a_{\tau}}\|^2-\frac{k\alpha}{\omega_2(1-\varphi_2)}\|e_a\|\|e_{a_{\tau}}\| \right), \\
        = &-\omega_2(1-\varphi_2)\left(\|e_{a_{\tau}}\|^2-\frac{k\alpha}{\omega_2(1-\varphi_2)}\|e_a\|\|e_{a_{\tau}}\| + \frac{k^2\alpha^2}{4\omega_2^2(1-\varphi_2)^2}\|e_a\|^2 \right) \\
        & +\frac{k^2\alpha^2}{4\omega_2(1-\varphi_2)}\|e_a\|^2, \\
        = &-\omega_2(1-\varphi_2)\left(\|e_{a_{\tau}}\|- \frac{k\alpha}{2\omega_2(1-\varphi_2)}\|e_a\| \right)^2 \\
        & +\frac{k^2\alpha^2}{4\omega_2(1-\varphi_2)}\|e_a\|^2, \\
        \leq & \frac{k^2\alpha^2}{4\omega_2(1-\varphi_2)}\|e_a\|^2.
    \end{split}
\end{equation} 
\normalsize
\noindent Second:
\begin{equation} \label{ineq 5}
    \begin{split}
        &k\|e_u\| \|e_{a_{\hat{\tau}}}\|-\omega_1\|e_{a_{\hat{\tau}}}\|^2 \\
        = & -\omega_1\left(\|e_{a_{\hat{\tau}}}\|^2 - \frac{k}{\omega_1}\|e_u\|\|e_{a_{\hat{\tau}}}\| \right),\\
        = & -\omega_1\left(\|e_{a_{\hat{\tau}}}\|^2 - \frac{k}{\omega_1}\|e_u\|\|e_{a_{\hat{\tau}}}\| +\frac{k^2}{4\omega_1^2}\|e_u\|^2 \right) + \frac{k^2}{4\omega_1}\|e_u\|^2, \\
        = & -\omega_1\left(\|e_{a_{\hat{\tau}}}\| - \frac{k}{2\omega_1}\|e_u\|^2 \right)^2 + \frac{k^2}{4\omega_1}\|e_u\|^2, \\ 
        \leq & \frac{k^2}{4\omega_1}\|e_u\|^2. \\
    \end{split}
\end{equation} 

Additionally, using Young's Inequality, we can obtain two more inequalities:
\begin{align} \label{ineq 6}
    &c_1\|e_a\| \leq \frac{1}{k\alpha}c_1^2+\frac{k\alpha}{4}\|e_a\|^2;\text{ and} \\
    &\|e_a\|\|z\|\rho(\|z\|) \leq \frac{1}{k\alpha}\rho^2(\|z\|)\|z\|^2+\frac{k\alpha}{4}\|e_a\|^2. \label{ineq 7}
\end{align}

By using the inequalities \eqref{ineq 4}--\eqref{ineq 7}, the expression in \eqref{vdot 4} can be upper bounded as 

\small
\begin{equation}\label{vdot 5}
    \begin{split}
        \dot{V} \leq &-\sum_{i=1}^{n-2}\|e_i\|^2 -\left(1-\frac{\varepsilon_2}{2}\right)\|e_{n-1}\|^2  -\left(\lambda-\left(\frac{\alpha}{2\varepsilon_1}+\frac{1}{2\varepsilon_2}\right) \right) \|e_n\|^2\\ &+\left(\frac{\alpha\varepsilon_1}{2}+\frac{k^2}{4\omega_1}\right)\|e_u\|^2 + \frac{\alpha}{2}\|\dot{u}_{\hat{\tau}} - \dot{u}_{\tau}\|^2 \\
        & -\bigg(\alpha k - \frac{\alpha}{2}-\frac{k^2\alpha^2}{4\omega_2\left(1-\varphi_2\right)} -\frac{k\alpha}{4}-\frac{k\alpha}{4}  \\
        &  -\omega_1 -\omega_2-\omega_3\left(\Bar{\Tilde{\tau}}+\hat{\tau}\right)k\bigg) \|e_a\|^2 \\ 
        & +\frac{1}{k\alpha}c_1^2 +\frac{1}{k\alpha}\rho^2(\|z\|)\|z\|^2 - \omega_3\int_{t-(\Bar{\Tilde{\tau}}+\hat{\tau})}^t\| \dot{u}(\theta) \|^2 d\theta .
    \end{split}
\end{equation}
\normalsize

We use the Cauchy?Schwarz inequality to develop the following upper bound
\begin{equation} \label{cauchy}
    \|e_u\|^2\leq \hat{\tau} \int_{t-\hat{\tau}}^t\|\dot{u}(\theta)\|^2d\theta. 
\end{equation}
By using \eqref{cauchy}, we have
\begin{equation} \label{ineq w1}
    \begin{split}
        \|e_u\|^2 &\leq \hat{\tau} \int_{t-(\Bar{\Tilde{\tau}}+\hat{\tau})}^t\|\dot{u}(\theta)\|^2d\theta, \\
        \Rightarrow\frac{-\omega_3}{4\hat{\tau}}\|e_u\|^2 &\geq -\frac{\omega_3}{4}\int_{t-(\Bar{\Tilde{\tau}}+\hat{\tau})}^t\|\dot{u}(\theta)\|^2d\theta.
    \end{split}
\end{equation}

Equation \eqref{udot} and $\Tilde{\tau}=\tau-\hat{\tau}$ yield the following: 
\begin{align} \label{ineq 8}
    \int_{t-\hat{\tau}}^t\|\dot{u}(\theta)\|^2d\theta \leq& k^2\int_{t-(\Bar{\Tilde{\tau}}+\hat{\tau})}^t\|e_a(\theta)\|^2, \\
    \int_{t-\tau}^t\|\dot{u}(\theta)\|^2d\theta \leq & k^2\int_{t-(\Bar{\Tilde{\tau}}+\hat{\tau})}^t\|e_a(\theta)\|^2, \label{ineq 9} \\
    \int_{t-(\Bar{\Tilde{\tau}}+\hat{\tau})}^t\|\dot{u}(\theta)\|^2d\theta =& k^2\int_{t-(\Bar{\Tilde{\tau}}+\hat{\tau})}^t\|e_a(\theta)\|^2d\theta. \label{ineq 10}
\end{align}

Using \eqref{Q1}, \eqref{ineq 8}, and \eqref{ineq 10}, we can write the following inequality 
\begin{equation} \label{ineq w2}
    \begin{split}
        \frac{Q_1}{\omega_1}k^2=k^2\int_{t-\hat{\tau}}^t\|e_a(\theta)\|^2d\theta &\leq k^2\int_{t-(\Bar{\Tilde{\tau}}+\hat{\tau})}^t\|e_a(\theta)\|^2, \\
        \frac{Q_1}{\omega_1}k^2&\leq \int_{t-(\Bar{\Tilde{\tau}}+\hat{\tau})}^t\|\dot{u}(\theta)\|^2d\theta, \\
        -\frac{\omega_3k^2}{4\omega_1}Q_1 &\geq -\frac{\omega_3}{4}\int_{t-(\Bar{\Tilde{\tau}}+\hat{\tau})}^t\|\dot{u}(\theta)\|^2d\theta.
    \end{split}
\end{equation}

Similarly using \eqref{Q2}, \eqref{ineq 9}, and \eqref{ineq 10}, we can write the following inequality

\begin{equation} \label{ineq w3}
    \begin{split}
        \frac{Q_2}{\omega_2}k^2=k^2\int_{t-\tau}^t\|e_a(\theta)\|^2d\theta &\leq k^2\int_{t-(\Bar{\Tilde{\tau}}+\hat{\tau})}^t\|e_a(\theta)\|^2, \\
        \frac{Q_2}{\omega_2}k^2&\leq \int_{t-(\Bar{\Tilde{\tau}}+\hat{\tau})}^t\|\dot{u}(\theta)\|^2d\theta, \\
        -\frac{\omega_3k^2}{4\omega_2}Q_2 &\geq -\frac{\omega_3}{4}\int_{t-(\Bar{\Tilde{\tau}}+\hat{\tau})}^t\|\dot{u}(\theta)\|^2d\theta.
    \end{split}
\end{equation}

\noindent From \eqref{Q3}, the following bound can be obtained for $Q_3$: 
\begin{equation} \label{ineq 11}
    \begin{split}
        Q_3 &\leq \omega_3(\Bar{\Tilde{\tau}}+\hat{\tau}) \sup_{s\in[t-(\Bar{\Tilde{\tau}}+\hat{\tau}),t]} \int_s^t\|\dot{u}(\theta)\|^2 d\theta, \\
        & \leq \omega_3(\Bar{\Tilde{\tau}}+\hat{\tau}) \int_{t-(\Bar{\Tilde{\tau}}+\hat{\tau})}^t\|\dot{u}(\theta)\|^2d\theta.
    \end{split}
\end{equation}
Using \eqref{ineq 11}, we obtain the following inequality: 
\begin{equation} \label{ineq w4}
    \begin{split}
        \frac{Q_3}{\omega_3(\Bar{\Tilde{\tau}}+\hat{\tau})} & \leq \int_{t-(\Bar{\Tilde{\tau}}+\hat{\tau})}^t\|\dot{u}(\theta)\|^2d\theta,\\
        \Rightarrow -\frac{1}{4(\Bar{\Tilde{\tau}}+\hat{\tau})}Q_3 &\geq -\frac{\omega_3}{4}\int_{t-(\Bar{\Tilde{\tau}}+\hat{\tau})}^t\|\dot{u}(\theta)\|^2d\theta.
    \end{split}
\end{equation}

Define $\zeta \triangleq \omega_1 + \omega_2 +\omega_3k(\Bar{\Tilde{\tau}}+\hat{\tau})+\frac{\alpha}{2}$. Substituting for $\omega_2$ and using the inequalities \eqref{ineq w1}, \eqref{ineq w2}, \eqref{ineq w3}, and \eqref{ineq w4}, the expression \eqref{vdot 5} can be upper bounded as

\small
\begin{equation} \label{vdot 6}
    \begin{split}
        \dot{V} \leq &-\sum_{i=1}^{n-2}\|e_i\|^2 -(1-\frac{\varepsilon_2}{2})\|e_{n-1}\|^2  -\left(\lambda-(\frac{\alpha}{2\varepsilon_1}+\frac{1}{2\varepsilon_2}) \right) \|e_n\|^2\\ &-\left(\frac{\omega_3}{4\hat{\tau}}-\left(\frac{\alpha\varepsilon_1}{2}+\frac{k^2}{4\omega_1}\right)\right)\|e_u\|^2 -\frac{k\alpha}{8}\|e_a\|^2 \\
        &  -\left(\frac{k\alpha}{8} - \zeta\right)\|e_a\|^2 + \frac{\alpha}{2}\|\dot{u}_{\hat{\tau}} - \dot{u}_{\tau}\|^2 +\frac{1}{k\alpha}c_1^2 +\frac{1}{k\alpha}\rho^2(\|z\|)\|z\|^2\\ 
        &  -\frac{\omega_3k^2}{4\omega_1}Q_1  -\frac{\omega_3k^2}{4\omega_2}Q_2 -\frac{1}{4(\Bar{\Tilde{\tau}}+\hat{\tau})}Q_3.
    \end{split}
\end{equation}
\normalsize
To show that $\dot{V}$ is negative semi-definite for the Lyapunov analysis, and given that $\|.\|^2\geq 0$, terms $(1-\frac{\varepsilon_2}{2})$, $\left(\lambda-(\frac{\alpha}{2\varepsilon_1}+\frac{1}{2\varepsilon_2}) \right)$, $\left(\frac{\omega_3}{4\hat{\tau}}-\left(\frac{\alpha\varepsilon_1}{2}+\frac{k^2}{4\omega_1}\right)\right)$, and $\left(\frac{k\alpha}{8} - \zeta\right)$ in \eqref{vdot 6} must be positive. This motivates the gain conditions in Theorem \ref{theorem 1}. 

Using Proposition 3 in \cite{cortes2008discontinuous}, there exists a Filippov solution \cite{filippov2013differential} for \eqref{udot}. Therefore, from \eqref{udot}, we have
\begin{align}
    \ddot{u}(t)&=\Big(\delta(e_1(t))\delta(\operatorname{sgn}(e_1(t))+1)\dot{e}_1(t)\Big)k e_a(t) \nonumber\\ 
    &~~~~+ \operatorname{sgn}\left(\frac{\operatorname{sgn}(e_1(t))+1}{2}\right) k \dot{e}_a(t), \\
&= \operatorname{sgn}\left(\frac{\operatorname{sgn}(e_1(t))+1}{2}\right) k \dot{e}_a(t). \label{u_ddot}
\end{align}
Next, using Assumption \ref{assumption error bound} and \eqref{e_a error}, it can be shown that all terms in \eqref{rdot3} are bounded. Hence, we have 
\begin{equation} \label{u_ddot bound}
    \ddot{u}(t) <m,
\end{equation} 
where $m\in\mathbb{R}_+$ is some positive constant.

Using \eqref{u_ddot bound} and the Mean Value Theorem, we can obtain $\|\dot{u}_{\hat{\tau}} - \dot{u}_{\tau}\|\leq \|\ddot{u}(\Theta(t,\hat{\tau}))\||\Tilde{\tau}| \leq m|\Tilde{\tau}|$, where $\Theta(t,\hat{\tau}))$ is between $t-\tau$ and $t-\hat{\tau}$, and $m$ is a positive constant.
Using \eqref{sigma bound}, \eqref{gain conditions}, and the inequality $\|z\|\leq\|y\|$, we obtain the following bound for \eqref{vdot 6}:

\begin{equation} \label{vdot 7}
    \begin{split}
        \dot{V} \leq & -\left(\frac{\sigma}{2}-\frac{1}{k\alpha}\rho^2(\|y\|) \right)\|z\|^2 - \frac{\sigma}{2}\|z\|^2 \\
        &  -\frac{\omega_3k^2}{4\omega_1}Q_1  -\frac{\omega_3k^2}{4\omega_2}Q_2 -\frac{1}{4(\Bar{\Tilde{\tau}}+\hat{\tau})}Q_3\\
        &+ \frac{\alpha \Bar{\Tilde{\tau}}^2m^2}{2} +\frac{1}{k\alpha}c_1^2.
    \end{split}
\end{equation}

\noindent By utilizing \eqref{delta bound}, the expression \eqref{vdot 7} reduces to

\begin{equation}
    \dot{V} \leq -\Delta\|y\|^2, \,\, \forall\|y\|\geq \sqrt{\frac{2c_1^2+k\alpha^2\Bar{\Tilde{\tau}}^2m^2}{2k\alpha\Delta}}. 
\end{equation}

\noindent Using \eqref{lyap ineq}, we conclude that $y$ is uniformly ultimately bounded, in the sense that $\limsup_{t\rightarrow\infty}\|y(t)\|\leq \sqrt{\frac{2c_1^2+k\alpha^2\Bar{\Tilde{\tau}}^2m^2}{k\alpha\Delta}}.$

\end{document}